\documentclass[final]{siamltex}
\usepackage{xcolor}
\usepackage{amsmath}
\usepackage{amssymb}
\usepackage{mathabx}
\usepackage{mathrsfs}
\usepackage{graphicx}
\usepackage{bm}
\usepackage{ucs}
\usepackage[utf8x]{inputenc}
\usepackage{wrapfig}
\usepackage{color}
\usepackage{float}
\usepackage[breaklinks,colorlinks=true,linkcolor=blue,citecolor=red,
 backref=page]{hyperref}

\newcommand{\ea}{\end{eqnarray}}  
\newcommand{\ba}{\begin{eqnarray}}  
\newcommand{\ee}{\end{equation}}  
\newcommand{\be}{\begin{equation}}  
\newcommand{\ean}{\end{eqnarray*}}  
\newcommand{\ban}{\begin{eqnarray*}}  
\renewcommand{\tilde}{\widetilde}

\DeclareMathAlphabet{\itbf}{OML}{cmm}{b}{it}
\newtheorem{remark}{Remark}[section]
\DeclareMathAlphabet{\itbf}{OML}{cmm}{b}{it}
 \DeclareMathAlphabet\mathbfcal{OMS}{cmsy}{b}{n}

\def\EE{\mathbb{E}}

\def\RR{\mathbb{R}}
\def\eps{\varepsilon}

\def\vr{{\vec{\itbf r}}}
\def\vq{{\vec{\itbf q}}}
\def\br{{{\itbf r}}}
\def\by{{{\itbf y}}}
\def\bv{{{\itbf v}}}
\def\bk{{{\itbf k}}}
\def\bK{{{\itbf K}}}
\def\bx{{\itbf x}}
\def\bX{{\itbf X}}
\def\bq{{{\itbf q}}}

\def\tae{A^\varepsilon}

\newcommand{\bnu}{{\boldsymbol{\nu}}}
\newcommand{\vnu}{{\vec{\bnu}}}
\newcommand{\bxi}{{\boldsymbol{\xi}}}

\newcommand{\la}{\lambda}

\newcommand{\vx}{\vec{\itbf x}}

\newcommand{\om}{\omega}

\newcommand{\cR}{{\mathscr R}}

\newcommand{\ep}{\varepsilon}

\newcommand{\vv}{\vec{\bv}}

\renewcommand{\hat}{\widehat}

\begin{document}

\title{Wave propagation and imaging in moving random media}

\author{
Liliana Borcea\footnotemark[1]
\and
Josselin Garnier\footnotemark[2]
\and 
Knut Solna\footnotemark[3]
}

\maketitle

\renewcommand{\thefootnote}{\fnsymbol{footnote}}

\footnotetext[1]{Department of Mathematics, University of Michigan,
  Ann Arbor, MI 48109. {\tt borcea@umich.edu}}
\footnotetext[2]{Centre de Math\'ematiques Appliqu\'ees, Ecole
  Polytechnique, 91128 Palaiseau Cedex, France.  {\tt
    josselin.garnier@polytechnique.edu}}
\footnotetext[3]{Department of Mathematics,
University of California at Irvine,
Irvine, CA 92697. {\tt ksolna@math.uci.edu}}

\begin{abstract}
We present a study of sound wave propagation in a time dependent
random medium and an application to imaging. The medium is modeled by
small temporal and spatial random fluctuations in the wave speed and
density, and it moves due to an ambient flow. We develop a transport
theory for the energy density of the waves, in a forward scattering
regime, within a cone (beam) of propagation with small opening
angle. We apply the transport theory to the inverse problem of
estimating a stationary wave source from measurements at a remote
array of receivers. The estimation requires knowledge of the mean
velocity of the ambient flow and the second-order statistics of the
random medium.  If these are not known, we show how they may be
estimated from additional measurements gathered at the array, using a
few known sources. We also show how the transport theory can be used
to estimate the mean velocity of the medium. If the array has large
aperture and the scattering in the random medium is strong, this
estimate does not depend on the knowledge of the statistics of the
random medium.
 \end{abstract}


\section{Introduction}
\begin{figure}[t]
\begin{picture}(0,170)%
\hspace{1.3in}\includegraphics[width = 0.5\textwidth]{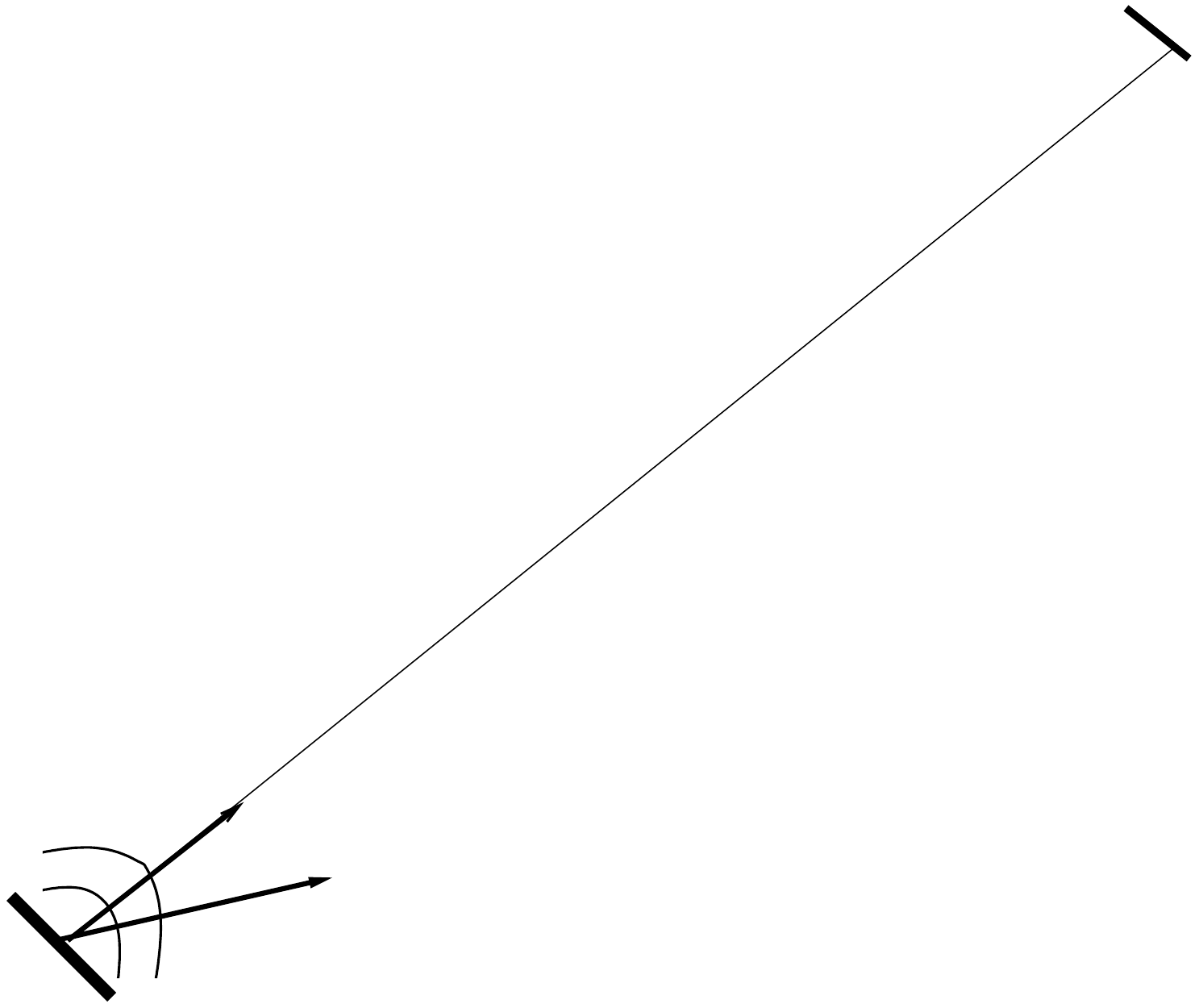}%
\end{picture}%
\setlength{\unitlength}{3947sp}%
\begingroup\makeatletter\ifx\SetFigFont\undefined%
\gdef\SetFigFont#1#2#3#4#5{%
  \reset@font\fontsize{#1}{#2pt}%
  \fontfamily{#3}\fontseries{#4}\fontshape{#5}%
  \selectfont}%
\fi\endgroup%
\begin{picture}(9770,1095)(1568,-1569)
\put(6050,1000){\makebox(0,0)[lb]{\smash{{\SetFigFont{7}{8.4}{\familydefault}{\mddefault}{\updefault}{\color[rgb]{0,0,0}{\normalsize $\mbox{receiver array}$}}%
}}}}
\put(3650,-1000){\makebox(0,0)[lb]{\smash{{\SetFigFont{7}{8.4}{\familydefault}{\mddefault}{\updefault}{\color[rgb]{0,0,0}{\normalsize $z$}}%
}}}}
\put(3850,-1180){\makebox(0,0)[lb]{\smash{{\SetFigFont{7}{8.4}{\familydefault}{\mddefault}{\updefault}{\color[rgb]{0,0,0}{\normalsize $\vv_o$}}%
}}}}
\put(2700,-1480){\makebox(0,0)[lb]{\smash{{\SetFigFont{7}{8.4}{\familydefault}{\mddefault}{\updefault}{\color[rgb]{0,0,0}{\normalsize $\mbox{source}$}}%
}}}}
\end{picture}%
\caption{Illustration of the setup. A stationary source emits a wave
  in the range direction $z$, in a moving medium with velocity
  $\vv(t,\vx)$ that has small random fluctuations about the constant
  mean $\vv_o$. The orientation of $\vv_o$ with respect to the range
  direction is arbitrary. The wave is recorded by a stationary, remote
  array of receivers. }
\label{fig:setup}
\end{figure}
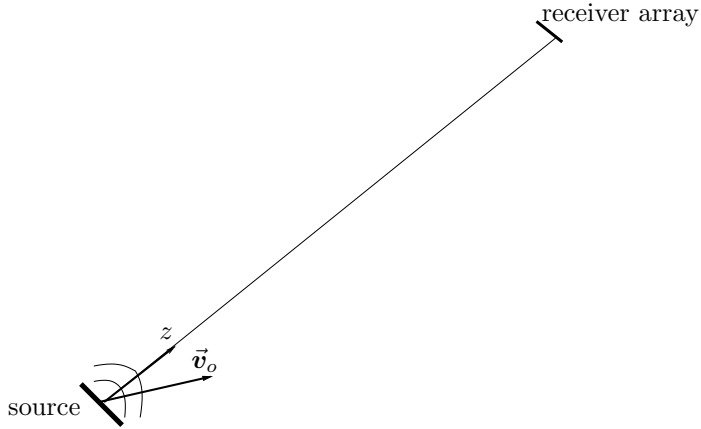

We study sound wave propagation in a time dependent medium modeled by
the wave speed $c(t,\vx)$ and density $\rho(t,\vx)$ that are random
perturbations of the constant values $c_o$ and $\rho_o$.  The medium
is moving due to an ambient flow, with velocity $\vv(t,\vx)$ that has
a constant mean $\vv_o$ and small random fluctuations. The source is
at a stationary location and emits a signal in the range direction
denoted henceforth by the coordinate $z$, as illustrated in Figure
\ref{fig:setup}. The signal is typically a pulse defined by an envelope
function of compact support, modulated at frequency $\om_o$.  It
generates a wave that undergoes scattering as it propagates through
the random medium. The goal of the paper is to analyze from first
principles the net scattering at long range, and to apply the results
to the inverse problem of estimating the source location and medium
velocity from measurements of the wave at a remote, stationary array
of receivers.

Various models of sound waves in moving media are described in
\cite[Chapter 2]{ostashev2015acoustics} using the linearization of the
fluid dynamics equations about an ambient flow, followed by
simplifications motivated by scaling assumptions. Here we consider
Pierce's equations \cite[Section 2.4.6]{ostashev2015acoustics} derived
in \cite{pierce1990wave} for media that vary at longer scales than the
central wavelength $\la_o = 2 \pi c_o/\om_o$ of the wave generated by
the source. Pierce's model gives the acoustic pressure
\begin{equation}
p(t,\vx) = - \rho(t,\vx) D_t \phi(t,\vx),  
\label{eq:P1}
\end{equation}
in terms of the velocity quasi-potential $\phi(t,\vx)$, which
satisfies the equation
\begin{equation}
D_t \Big[ \frac{1}{c^2(t,\vx)} D_t \phi(t,\vx)\Big] -
\frac{1}{\rho(t,\vx)} \nabla_{\vx} \cdot \Big[ \rho(t,\vx) \nabla_{\vx}
  \phi(t,\vx)\Big] = {\rm s}(t,\vx),
\label{eq:P2}
\end{equation}
for spatial variable $\vx = (\bx,z) \in \RR^{d+1}$ and time $t \in
\RR$, with natural number $d \ge 1$.  Here $\bx \in \RR^d$ lies in the
cross-range plane, orthogonal to the range axis $z$. Moreover,
$\nabla_{\vx}$ and $\nabla_{\vx} \cdot$ are the gradient and
divergence operators in the variable $\vx$ and
\[
D_t = \partial_t + \vv(t,\vx) \cdot \nabla_{\vx}
\] 
is the material (Lagrangian) derivative, with $\partial_t$ denoting
the partial derivative with respect to time. The source is modeled by
the function ${\rm s}(t,\vx)$ localized at the origin of range and with
compact support. Prior to the source excitation there is no wave
\begin{equation}
\phi(t,\vx) \equiv 0, \quad t \ll 0,
\label{eq:P3}
\end{equation}
but the medium is in motion due to the ambient flow. 

Sound wave propagation in ambient flows due to wind in the atmosphere
or ocean currents arises in applications like the quantification of the
effects of temperature fluctuations and wind on the rise time and
shape of sonic booms \cite{boulanger1995sonic} or on radio-acoustic
sounding \cite{karyukin1982influence}, monitoring noise near airports
\cite{schomer1983noise}, acoustic tomography \cite{munk2009ocean}, and
so on.

Moving media also arise in optics, for example in Doppler velocimetry
or anemometry \cite{drain1980laser, durstprinciples} which uses lasers
to determine the flow velocity $\vv_o$. This has applications in wind
tunnel experiments for testing aircraft \cite{garman2006airborne}, in
velocity analysis of water flow for ship hull design
\cite{molland2017ship}, in navigation and landing
\cite{amzajerdian2013lidar}, in medicine and bioengineering
\cite{oberg1990laser}. A description of light propagation models used
in this context can be found in \cite[Chapter 8]{Mayinger}.

Much of the applied literature on waves in moving random media
considers either discrete models with Rayleigh or Mie scattering by
moving particles \cite{Mayinger} or continuum models described by the
classic wave equation with wave speed $c(0,\vx-\vec{\bv}t)$. These use
Taylor's hypothesis \cite[Chapter 19]{akira1978wave} where the medium
is ``frozen" over the duration of the experiment {and simply shifted by the uniform ambient flow}.  
A transport theory in such frozen-in media is obtained for example in \cite[Chapter
  20]{akira1978wave} and \cite[Chapter 8]{ostashev2015acoustics}, in
the paraxial regime where the waves propagate in a narrow angle cone
around the range direction. The formal derivation of this theory
assumes that the random fluctuations of the wave speed are Gaussian,
and uses the Markov approximation, where the fluctuations are
$\delta$--correlated in range i.e, at any two distinct ranges, no
matter how close, the fluctuations are assumed uncorrelated.

In this paper we study the wave equation \eqref{eq:P2} with
coefficients $c(t,\vx), \rho(t,\vx)$ and $\vv(t,\vx)$ that have random
correlated fluctuations at spatial scale $\ell$ and temporal scale
$T$. These fluctuations are not necessarily Gaussian. We analyze the
solution $\phi(t,\vx)$ and therefore the acoustic pressure $p(t,\vx)$
in a forward scattering regime, where the propagation is within a cone
(beam) with axis along the range direction $z$. The analysis uses
asymptotics in the small parameter $\ep = \la_o/L \ll 1$, where $L$ is
the range scale that quantifies the distance between the source and
the array of receivers.  Pierce's equations
\eqref{eq:P1}--\eqref{eq:P2} are justified for small $\la_o/\ell \ll
1$.  By fixing $\la_o/\ell$ or letting it tend to zero, independent of
$\ep$, and by appropriate scaling of the spatial support of the
source ${\rm s}(t,\vx)$, we obtain two wave propagation
regimes: The first is called the wide beam regime because the cone of
propagation has finite opening angle. The second is the paraxial
regime, where the cone has very small opening angle.  We use the
diffusion approximation theory given in \cite[Chapter 6]{book} and
\cite{PK1974,PW1994} to study both regimes and obtain transport
equations that describe the propagation of energy. These equations are
simpler in the paraxial case and we use them to study the inverse
problem of locating the source. Because the inversion requires
knowledge of the mean velocity $\vv_o$ of the ambient flow and the
second-order statistics of the random medium, we also discuss their
estimation from additional measurements of waves generated by known
sources.

The paper is organized as follows: We begin in section \ref{sect:form}
with the mathematical formulation of the problem. Then we state in
section \ref{sect:ResForward} the transport equations. These equations
are derived in section \ref{sect:Anal} and we use them for the inverse
problem in section \ref{sect:Inverse}. We end with a summary in
section \ref{sect:sum}.

\section{Formulation of the problem}
\label{sect:form}
We study the sound wave modeled by the acoustic pressure $p(t,\vx)$
defined in equation \eqref{eq:P1} in terms of the velocity
quasi-potential $\phi(t,\vx)$, the solution of the initial value
problem \eqref{eq:P2}--\eqref{eq:P3}. 
The problem is  to characterize the acoustic pressure $p(t,\vx)$
in the scaling regime described in section \ref{sect:scaling} and then use the results for
localizing the source and estimating the mean medium velocity $\vv_o$.

\subsection{Medium and source}
The coefficients in equation
\eqref{eq:P2} are random fields, defined by 
\begin{align}
\vv(t,\vx) &= \vv_o + V \sigma_v \vnu
\Big(\frac{t}{T},\frac{\vx-\vv_o t}{\ell} \Big), \label{eq:F1}
\\ \rho(t,\vx) &= \rho_o \exp \Big[ \sigma_\rho
  \nu_\rho\Big(\frac{t}{T},\frac{\vx-\vv_o t}{\ell} \Big)
  \Big], \label{eq:F2} \\ c(t,\vx) &= c_o \Big[ 1 + \sigma_c
  \nu_c\Big(\frac{t}{T},\frac{\vx-\vv_o t}{\ell} \Big) \Big]^{-1/2},
 \label{eq:F3}
\end{align}
where $c_o$, $\rho_o$ are the constant background wave speed and
density, $\vv_o$ is the constant mean velocity of the ambient flow, and
$V$ is a velocity scale  (of the order of $|\vv_o|$) that will be specified later. The fluctuations in
\eqref{eq:F1}--\eqref{eq:F3} are given by the random stationary
processes $\vnu$, $\nu_\rho$ and $\nu_c$ of dimensionless arguments
and mean zero
\begin{equation}
\EE\big[ \vnu(\tau,\vr)\big] = 0, \quad \EE[\nu_\rho(\tau,\vr)] = 0, \quad 
 \EE[\nu_c(\tau,\vr)] = 0.
\label{eq:F4}
\end{equation} We assume that $\vnu = (\nu_j)_{j=1}^{d+1}$, $\nu_\rho$ and
$\nu_c$ are twice differentiable, with bounded derivatives almost
surely, have ergodic properties in the $z$ direction, and are
correlated, with covariance entries
\begin{equation}
\EE\big[ \nu_\alpha (\tau, \vr) \nu_\beta (\tau',\vr') \big] =
\cR_{\alpha \beta} (\tau-\tau',\vr-\vr').
\label{eq:F5}
\end{equation}
Here the indices $\alpha$ and $\beta$ are either $1 ,\ldots, d+1$, or $\rho$, or $c$. The covariance is an even and integrable symmetric
matrix valued function, which is four times
differentiable
  and satisfies the
normalization conditions
\begin{equation}
\cR_{\alpha \alpha}(0,0) = 1 \mbox{ or } O(1), \quad \quad \int_{\RR} {\rm d} \tau
\int_{\RR^{d+1}} {\rm d} \vr \, \cR_{\alpha \alpha}(\tau,\vr) = 1 ~
\mbox{ or } O(1).
\label{eq:F6}
\end{equation}
The scale $T$ in definitions \eqref{eq:F1}--\eqref{eq:F3} is the
correlation time, the typical lifespan of a spatial realization of the
fluctuations, and $\ell$ is the correlation length, the typical length
scale of the fluctuations.  The dimensionless positive numbers
$\sigma_v$, $\sigma_\rho$ and $\sigma_c$ quantify the standard
deviation of the fluctuations. They are of the same order and small,
so definitions \eqref{eq:F2} and \eqref{eq:F3} can be approximated by
\begin{align*}
\rho(t,\vx) & \approx \rho_o \Big[ 1 + \sigma_\rho
  \nu_\rho\Big(\frac{t}{T},\frac{\vx-\vv_o t}{\ell} \Big) \Big],
\quad  c(t,\vx) \approx c_o \Big[ 1 - \frac{\sigma_c}{2}
  \nu_c\Big(\frac{t}{T},\frac{\vx-\vv_o t}{\ell} \Big) \Big],
\end{align*}
with $c_o$ and $\rho_o$ close to the mean wave speed and density.  The
exponential in \eqref{eq:F2} and the inverse of the square root in
\eqref{eq:F3} are used for convenience because some important effective properties 
of the medium are defined in terms of $\EE[\log \rho]$ and $\EE[c^{-2}]$,
which are equal to $\log \rho_o$ and $c_o^{-2}$.

The origin of the coordinates is at the center of the 
source location, modeled by
\begin{equation}
{\rm s}(t,\vx) = \sigma_s e^{-i \om_o t} S
\Big(\frac{t}{T_s},\frac{\bx}{\ell_s} \Big) \delta(z),
\label{eq:F7}
\end{equation}
for $\vx = (\bx,z)$, using the continuous function $S$ of
dimensionless arguments and compact support. The length scale $\ell_s$
is the radius of the support of ${\rm s}(t,\vx)$ in cross-range and the time
scale $T_s$ is the duration of the emitted signal.  Note that ${\rm
  s}(t,\vx)$ is modulated by the oscillatory exponential at the
frequency $\om_o$. We call it the central frequency because the
Fourier transform of ${\rm s}(t,\vx)$ with respect to time is
supported in the frequency interval $|\om-\om_o| \le O(1/T_s)$. The
solution $\phi(t,\vx)$ of \eqref{eq:F2} depends linearly on the
source, so we use $\sigma_s$ to control its amplitude.

To be able to set radiation conditions for  the wave field resolved over
frequencies, we make the mathematical assumption that the random
fluctuations of $\vv(t,\vx)$, $\rho(t,\vx)$ and $c(t,\vx)$ are
supported in a domain of finite range that is much larger than $L$. In
practice this assumption does not hold, but the wave equation is
causal and with finite speed of propagation, so the truncation of the
support of the fluctuations does not affect the wave measured at the
array up to time $O(L/c_o)$.

\subsection{Scaling regime}
\label{sect:scaling}
Because the fluctuations of the coefficients
\eqref{eq:F1}--\eqref{eq:F2} are small, they have negligible effect on
the wave at short range, meaning that $\phi(t,\vx) \approx
\phi_o(t,\vx)$, the solution of \eqref{eq:P2}--\eqref{eq:P3} with
constant wave speed $c_o$, density $\rho_o$ and velocity $\vv_o$. We
are interested in a long range $L$, where the wave undergoes many
scattering events in the random medium and $\phi(t,\vx)$ is quite
different from $\phi_o(t,\vx)$.  We model this long range regime with
the small and positive, dimensionless parameter
\begin{equation}
\ep = \frac{\la_o}{L} \ll 1,
\label{eq:F8}
\end{equation}
and use asymptotics in the limit $\ep \to 0$ to study the random field 
$\phi(t,\vx)$.

The relation between the wavelength, the correlation length and the
cross-range support of the source is described by the positive,
dimensionless parameters
\begin{equation}
\gamma = \frac{\la_o}{\ell}, \quad \gamma_s = \frac{\la_o}{\ell_s},
\label{eq:F9}
\end{equation}
which are small, but independent of $\ep$. The
positive, dimensionless parameter
\begin{equation}
\eta = \frac{T}{T_L}, 
\label{eq:F10}
\end{equation}
determines how fast the medium changes on the scale of the travel time
$T_L = L/c_o$.

The duration of the source signal is modeled by the positive,
dimensionless parameter
\begin{equation}
\eta_s = \frac{T_s}{T_L},
\label{eq:F11}
\end{equation}
which is independent of $\ep$. The Fourier transform of this signal is
supported in the frequency interval centered at $\om_o$ and of length
(bandwidth) $O(1/T_s)$, where
\begin{equation}\label{eq:b2}
 \frac{1}{T_s} = \frac{1}{\eta_s T_L}  \ll  \frac{1}{\ep T_L} 
   = \frac{c_o}{\ep L}  = \frac{c_o}{\lambda_o} = O(\omega_o)  .
\end{equation}
Thus, the source has  a small bandwidth in the $\ep \to 0$ limit.

Our asymptotic analysis assumes the order relation 
\begin{equation}
\ep \ll \min \{\gamma, \gamma_s, \eta, \eta_s\}   ,
\end{equation}
meaning that we take the limit $\ep \to 0$ for fixed $\gamma,
\gamma_s, \eta, \eta_s$. The standard deviations of the fluctuations
are scaled as
\begin{equation}
\sigma_c = \sqrt{\ep \gamma} \bar\sigma_c, \quad\sigma_\rho = \sqrt{\ep \gamma} \bar\sigma_\rho, \quad
\sigma_v = \sqrt{\ep \gamma} \bar\sigma_v,
\label{eq:F12}
\end{equation}
with $\bar\sigma_c,\bar\sigma_\rho,\bar\sigma_v=O(1)$
to obtain a $O(1)$ net scattering effect.

The ambient flow, due for example to wind, has much smaller velocity
than the reference sound speed $c_o$. We model this assumption with the
scaling relation
\begin{equation}
|\vv_o|/V = O(1), \quad \mbox{where } V = \ep c_o.
\label{eq:F13}
\end{equation}
Although $V \ll c_o$, the medium moves on the scale
of the wavelength over the duration of the propagation
\begin{equation}\label{eq:b1}
    V T_L = V \frac{L}{c_o} = \ep L =
    \lambda_o <  \ell,
\end{equation}
so the motion has a $O(1)$ net scattering effect.  Slower motion is
negligible, whereas faster motion gives different phenomena than those
analyzed in this paper.

We scale the amplitude of the source as 
\begin{equation}
\sigma_s = \frac{1}{\ep \eta_s  L} \Big(\frac{\gamma_s
}{\ep}\Big)^{d},
\label{eq:F14}
\end{equation}
to obtain $\phi(t,\vx) = O(1)$ in the limit $\ep \to 0$. Since
equation \eqref{eq:F2} is linear, any other source amplitude can be
taken into account by multiplication of our wave field with that given
amplitude.

Note that in section \ref{sect:resPar} we consider the secondary
scaling relation 
\begin{equation}
\gamma \sim \gamma_S \ll 1,
\label{eq:scParax1}
\end{equation}
corresponding to the paraxial regime, where the symbol ``$\sim$''
means of the same order. Moreover, in section
\ref{sect:Inv2} we assume $\eta/\eta_s \ll 1$ corresponding to a
regime of statistical stability.  In this secondary scaling regime we
let \ba |\vv_o| = O\left(\frac{\ep c_o}{\eta \gamma}\right),
\ea to obtain the distinguished limit in which the medium velocity
impacts the quantities of interest.

\section{Results of the analysis of the wave field}
\label{sect:ResForward}
We show in section \ref{sect:Anal} that in the scaling regime
described in equations \eqref{eq:F8}--\eqref{eq:F14}, the 
pressure is given by
\begin{align}
p(t,\vx) \approx i \om_o {\rho_o} \int_{\mathcal{O}} \frac{{\rm d}\om
  {\rm d} \bk}{(2 \pi)^{d+1}} \frac{ a(\om,\bk,z) }{\sqrt{\beta(\bk)}}
e^{-i (\om_o + \om)t +i \vec{\bk} \cdot \vx},
\label{eq:B3}
\end{align}
for $\vx = (\bx,z)$ and $\mathcal{O} = \{ \om \in \RR \} \times \{ \bk
\in \RR^d, ~ |\bk| < k_o\}, $ where the approximation error vanishes
in the limit $\ep \to 0$. This expression is a Fourier synthesis of
forward propagating time-harmonic plane waves (modes) at frequency
$\om_o + \om$, with wave vectors $\vec{\bk}$ defined by
\begin{equation}
\vec{\bk} = \big(\bk,\beta(\bk)\big), \quad \beta(\bk) = \sqrt{k_o^2 -
  |\bk|^2}, \quad k_o = 2 \pi/\la_o.
\label{eq:R1}
\end{equation}
The scattering effects in the random medium
are captured by the mode amplitudes, which form a Markov process
$\big(a(\om,\bk,z)\big)_{(\om,\bk)\in \mathcal{O}}$ that evolves in
$z$, starting from
\begin{equation}
a(\om,\bk,0) = a_o(\om,\bk) = \frac{i \sigma_s T_s\ell_s^d }{2
  \sqrt{\beta(\bk)}} \hat S(\om T_s, \ell_s \bk).
\label{eq:R2}
\end{equation}
This process satisfies  the conservation relation
\begin{equation}
\int_{\mathcal{O}} {\rm d}\om {\rm d} \bk\, \big|a(\om,\bk,z)\big|^2 =
\int_{\mathcal{O}} {\rm d}\om {\rm d} \bk\, \big|a_o(\om,\bk)\big|^2,
\quad \forall z > 0.
\label{eq:R3}
\end{equation}

The statistical moments of 
$\big(a(\om,\bk,z)\big)_{(\om,\bk)\in \mathcal{O}}$ are characterized
explicitly in the limit $\ep \to 0$, as explained in section
\ref{sect:MarkLim} and Appendix \ref{ap:DL}. Here we describe the
expectation of the amplitudes, which defines the coherent wave, and
the second moments that define the mean Wigner transform of the wave
i.e., the energy resolved over frequencies and direction of
propagation.

\subsection{The coherent wave}
\label{sect:Res1}
The expectation of the acoustic pressure (the coherent wave) is
obtained from \eqref{eq:B3} using the mean amplitudes 
\begin{equation}
\EE[ a(\om,\bk,z)] = a_o(\om,\bk) \exp \left[ i \theta(\om,\bk) z +
  D( \bk) z\right].
\label{eq:R4}
\end{equation} 
These are derived in Section \ref{sect:mean}, with $a_o(\om,\bk)$
given in \eqref{eq:R2}. The exponential describes the effect of
the random medium, as follows:

The first term in the exponent is the phase 
\begin{equation}
\theta(\omega,\bk) = \frac{k_o}{\beta(\bk)} \Big(\frac{\om}{c_o} -
\frac{\bv_o}{c_o} \cdot \bk\Big) + 
\frac{ \sigma_\rho^2}{8 \beta(\bk) \ell^2}  
   \Delta_{\vr} \cR_{\rho \rho}(0,\vr)|_{\vr =
  0},
\label{eq:R5}
\end{equation}
and consists of two parts: The first part models the Doppler frequency
shift and depends on the cross-range component $\bv_o$ of the mean
velocity $\vec{\bv}_o=(\bv_o, v_{oz})$.  It comes from the
expansion of the mode wavenumber
\[
\sqrt{\Big(k_o + \frac{\om- \bv_o \cdot \bk}{c_o}\Big)^2 - |\bk|^2}
\approx \beta(\bk) + \frac{k_o (\om- \bv_o \cdot \bk)}{c_o \beta(\bk)}
,
\]
in the limit $\ep \to 0$, using the scaling relation (\ref{eq:F13})
and $\om \ll \om_o$ obtained from \eqref{eq:b2}. The second part is
due to the random medium and it is small when $\gamma \ll 1$, i.e., 
$\la_o \ll \ell$.

The second term in the exponent in \eqref{eq:R4} is  
\begin{align}
D(\bk) = - \frac{k_o^4  \ell^{d+1}}{4} \int_{|\bk'| < k_o}
\frac{{\rm d}\bk'}{(2\pi)^d} \, \frac{1}{\beta(\bk) \beta(\bk')}
\int_{\RR^d} {\rm d}\br \int_0^\infty {\rm d}r_z \, e^{- i \ell
  (\vec{\bk}-\vec{\bk}') \cdot \vr} \nonumber \\
   \times \Big[
\sigma_c^2  \cR_{cc}(0,\vr) + \frac{\sigma_{\rho}^2}{4 (k_o \ell)^4}
  \Delta_{\vr}^2 \cR_{\rho \rho}(0,\vr) - \frac{\sigma_{\rho}\sigma_c}{(k_o
    \ell)^2} \Delta_{\vr} \cR_{c \rho}(0,\vr) \Big], \label{eq:R6}
\end{align}
where we used the notation $\vr = (\br,r_z)$ and definition
\eqref{eq:R1}.  This complex exponent accounts for the significant
effect of the random medium, seen especially in the term proportional
to $\cR_{cc}$ which dominates the other ones in the $\gamma \ll 1$
regime. Because the covariance is even, the real part of $D(\bk)$
derives from
\begin{equation}
\int_{\RR^{d+1}} {\rm d} \vr \, \cR_{cc}(0,\vr) e^{-i
  \ell(\vec{\bk}-\vec{\bk'}) \cdot \vr} = \int_{\RR} \frac{{\rm d}
  \Omega}{2 \pi} \,
\widetilde\cR_{cc}\big(\Omega,\ell(\vec{\bk}-\vec{\bk'})\big),
\label{eq:TR4p}
\end{equation}
where 
\begin{equation}
\widetilde \cR_{cc}(\Omega,\vq) = \int_{\RR} {\rm d} \tau
\int_{\RR^{d+1}} {\rm d} \vr \, \cR_{cc}(\tau,\vr) e^{i \Omega \tau -
  i \vq \cdot \vr} \ge 0,
\label{eq:TR4}
\end{equation}
is the power spectral density of $\nu_c$. This is non-negative by
Bochner's theorem, so $\mbox{Re} \big[ D(\bk)\big] < 0$ and the mean
amplitudes decay exponentially in $z$, on the length scale
\begin{equation}
\mathscr{S}(\bk) = - \frac{1}{\mbox{Re} \big[ D(\bk)\big]},
\label{eq:R8}
\end{equation}
called the scattering mean free path.   Note that $|\bk|, |\bk'| =
O(1/\ell)$ in the support of $\widetilde \cR_{cc}$ in \eqref{eq:TR4p}
and that by choosing the standard deviation $\sigma_c$ as in
\eqref{eq:F12}, we obtain from \eqref{eq:R6}--\eqref{eq:R8} that
$\mathscr{S}(\bk) = O(L)$ in the $\ep \to 0$ followed by the $\gamma
\to 0$ limit.  This shows that the decay of the mean amplitudes in $z$
is significant in our regime. It is the manifestation of the
randomization of the wave due to scattering in the medium.
\subsection{The Wigner transform}
\label{sect:res2}
The strength of the random fluctuations of the mode amplitudes is
described by the Wigner transform (energy density)
\begin{align}
W(\omega,\bk,\bx,z) = \int \frac{{\rm d} {\itbf q}}{(2 \pi)^d} \, e^{i
  \itbf q \cdot (\nabla \beta(\bk) z + \bx)}
  \EE \Big[ a \Big(\omega,\bk + \frac{{\itbf q}}{2},z\Big) \overline{
      a\Big(\omega,\bk-\frac{ {\itbf q}}{2},z\Big)}\Big],
 \label{eq:T6:0}
\end{align}
where the bar denotes complex conjugate and the integral is over all
$\bq \in \RR^d$ such that $|\bk \pm \bq/2|<k_o$. The Wigner transform
satisfies the equation
\begin{align}
\big[ \partial_z - \nabla \beta(\bk) \cdot \nabla_\bx \big] W
(\omega,\bk,\bx,z) = \int_{\mathcal{O}}\frac{{\rm d}\omega'{\rm d}
  \bk'}{(2 \pi)^{d+1}}\, {\mathcal Q}(\omega,\omega',\bk,\bk') \big[
  W(\omega',\bk',\bx,z) \nonumber \\ -W(\omega,\bk,\bx,z) \big] ,
\label{eq:TRresult}
\end{align}
for $z>0$, with initial condition
\begin{equation}
\label{eq:TRinitial}
W (\omega,\bk,\bx,0) = |a_o(\om,\bk)|^2 \delta(\bx).
\end{equation}
The integral kernel in \eqref{eq:TRresult} is called the differential
scattering cross-section. It is defined by
\begin{align}
&{\mathcal Q}(\omega,\omega',\bk,\bk') = \frac{k_o^4    \ell^{d+1} T}{4 \beta(\bk)\beta(\bk') } \Big[
\sigma_c^2  \widetilde
\cR_{cc} + \frac{\sigma_{\rho}^2}{4(k_o \ell)^4}
\widetilde{ \Delta_{\scriptsize{\overrightarrow{\br}}}^2  {\cR}_{\rho
      \rho}}- \frac{\sigma_c \sigma_{\rho}}{(k_o \ell)^2}
 \widetilde{ \Delta_{\scriptsize{\overrightarrow{\br}}}  {\cR}_{c \rho} }
  \Big],
  \label{eq:DCS}
\end{align}
where the power spectral densities in the square bracket are evaluated
as \ba\label{eq:DCS2} \widetilde \cR_{cc} = \widetilde \cR_{cc}\big(T
(\om-\om' - (\vec{\bk}-\vec{\bk}') \cdot \vv_o), \ell
(\vec{\bk}-\vec{\bk'})\big), \ea and similar for the other two terms,
which are proportional to the Fourier transform of $\Delta_{\vr}^2
\cR_{\rho \rho}$ and $\Delta_{\vr} \cR_{c \rho}$.  The total
scattering cross section is defined by the integral of \eqref{eq:DCS}
and satisfies
\begin{align}
\Sigma(\bk) = \int_{\mathcal{O}}\frac{{\rm d} \omega'{\rm d} \bk'}{(2
  \pi)^{d+1}}\, {\mathcal Q}(\omega,\omega',\bk,\bk') =
  \frac{2}{{\mathscr{S}(\bk)}} .
\end{align}

Note that the last two terms in the square bracket in \eqref{eq:DCS}
are small in the $\gamma \ll 1$ regime, because $1/(k_o \ell) =
\gamma/(2 \pi) \ll 1$ and $ \sigma_{\rho}  / \sigma_c = O(1)$. If
$\sigma_{\rho} / \sigma_c$ were large, of the order $\gamma^{-2}$, then these
terms would contribute. However, this would only change the
interpretation of the differential scattering cross section and not
its qualitative form.

\subsubsection{The radiative transfer equation}
The evolution equation \eqref{eq:TRresult} for the Wigner transform is
related to the radiative transfer equation \cite{chandra,ryzhik96}.
Indeed, we show in Appendix \ref{app:tranproof} that
$W(\omega,\bk,\bx,z)$ is the solution of
\eqref{eq:TRresult}-\eqref{eq:TRinitial} if and only if
\begin{equation} \label{eq:T}
V(\omega,\vec\bk,\vec\bx) = \frac{1}{\beta(\bk)} W(\omega,\bk,\bx,z)
\delta\big( k_z - \beta(\bk)\big), \quad \vec \bk = (\bk,k_z),
\end{equation}
solves the radiative transfer equation
\begin{align}
\nonumber
\nabla_{\vec\bk} \Omega(\vec\bk) \cdot \nabla_{\vec\bx}
V(\omega,\vec\bk,\vec\bx) = \int_{\RR^{d+1}} \frac{{\rm d}
  \vec\bk'}{(2\pi)^{d+1}} \int\frac{{\rm d}\omega'}{2\pi}
\mathfrak{S}\big(\omega,\omega',\vec\bk,\vec\bk'\big) \big[
  V(\omega',\vec\bk',\vec\bx) \\- V(\omega,\vec\bk,\vec\bx) \big] ,
\label{eq:rte3d}
\end{align}
with $\Omega(\vec\bk)=c_o |\vec\bk|$ and  the scattering kernel
\begin{equation}
\mathfrak{S}\big(\omega,\omega',\vec\bk,\vec\bk'\big) = \frac{2 \pi
  c_o^2}{k_o^2} \beta(\bk) \beta(\bk') \mathcal{Q}(\om,\om',\bk,\bk')
\delta\big( \Omega(\vec\bk)-\Omega(\vec\bk')\big).
\end{equation}
The initial condition is specified at $\vx = (\bx,0)$ by
\begin{equation}
V(\om,\vec\bk,(\bx,0)) = |a_o(\om,\bk)|^2 \delta(\bx)
\delta(k_z-\beta(\bk)) ,
\label{eq:trIni}
\end{equation}
with $a_o(\om,\bk)$ defined in \eqref{eq:R2}.

This result shows that the generalized (singular) phase space energy \eqref{eq:T}
evolves as in the standard 3D radiative transfer equation, but it is
supported on the phase vectors with range component
\begin{equation}\label{eq:hs} 
k_z=\beta(\bk), \quad |\bk| <  k_o. 
\end{equation}
Indeed, if $\vec\bk' = (\bk',\beta(\bk'))$ and $\vec \bk= (\bk,k_z)$,
then
$$
\delta \big( \Omega(\vec\bk)- \Omega(\vec\bk') \big) = \frac{1}{c_o}
\delta \big(|\vec \bk|- k_o\big) = \frac{k_o}{c_o \beta(\bk)}
\delta\big(k_z -\beta(\bk)\big) ,
$$ so the evolution of $V(\om,\vec\bk,\vx)$ is confined to the
hypersurface in equation (\ref{eq:hs}). Physically, this means that
the wave energy is traveling with constant speed in a cone of
directions centered at the range axis $z$.

\subsubsection{Paraxial approximation}
\label{sect:resPar}
The paraxial approximation of the Wigner transform is obtained from
\eqref{eq:TRresult}-\eqref{eq:TRinitial} in the limit
\[
\gamma = \la_o/\ell \to 0, \quad \mbox{so that} ~ 
\gamma/\gamma_s = {\rm finite},
\]
as explained in section \ref{sect:Anal7}.  In this case the phase
space decomposition of the initial wave energy given by \eqref{eq:R2}
and \eqref{eq:TRinitial} is supported in a narrow cone around the
range axis $z$, with opening angle scaling as \ban
\frac{\la_o}{\ell_s} = \gamma_s \ll 1.  \ean Moreover, from the
expression \eqref{eq:DCS} of the differential scattering cross-section
and \eqref{eq:DCS2} we see that the energy coupling takes place in a
small cone of differential directions whose opening angle is \ban
\frac{\la_o}{\ell} = \gamma \ll 1.  \ean

In the paraxial regime equation \eqref{eq:TRresult} simplifies to
\begin{align}
\left[\partial_{z} + \frac{\bk}{k_o} \cdot \nabla_{\bx} \right]
  W(\omega,\bk,\bx,{z}) = \int_{\mathbb{R}^d} \frac{{\rm d} \bk'}{(2
    \pi)^d} \int_{\mathbb{R}}\frac{{\rm d} \omega'}{2 \pi} \,{\mathcal
    Q}_{\rm par}(\omega',\bk') \nonumber \\ \times W\left( \omega -
  \omega' - \bk' \cdot  \bv_o , \bk-\bk',\bx,{z}\right) -\Sigma_{{\rm
      par}} W(\omega,\bk,\bx,{z}), \label{eq:TRresultP}
\end{align}
where we obtained from definition \eqref{eq:R1} and the scaling
relations \eqref{eq:F10}, \eqref{eq:F13} that in the limit $\gamma \to
0$, \ban \beta(\bk) \to k_o, \quad \ell \big| \beta(\bk) -\beta(\bk')
\big| \to 0 , \quad T \big| v_{oz} ( \beta(\bk) -\beta(\bk') )\big|
\to 0.  \ean The differential scattering cross-section becomes
\begin{align}
{\mathcal Q}_{\rm par}(\omega, \bk ) = \frac{k_o^2 \sigma_c^2
  \ell^{d+1} T}{4 } \widetilde \cR_{cc} \big(T \omega , \ell\bk,
0\big), \label{eq:TrP2}
\end{align}
and the total scattering cross-section is
\begin{align}
\Sigma_{{\rm par}} = \int_{\mathbb{R}^d} \frac{{\rm d} \bk'}{(2
  \pi)^d} \int_{\mathbb{R}}\frac{{\rm d} \omega'}{2 \pi} \,{\mathcal
  Q}_{\rm par}(\omega',\bk') = \frac{ \sigma_c^2 \ell k_o^2
   }{4 } \cR(0,{\bf 0})=  \frac{2}{\mathscr{S}_{\rm par}}
,\label{eq:TrP3}
\end{align}
where $\mathscr{S}_{\rm par}$ is the scattering mean free path in the
paraxial regime and
\begin{equation}
\cR(\tau,\br) = \int_{\RR} {\rm d} r_z \, \cR_{cc}(\tau,\vr), \quad \vr = (\br,r_z).
\label{eq:Sol01}
\end{equation}
The initial condition is as in \eqref{eq:TRinitial}, with $a_o$
defined in \eqref{eq:R2},
\begin{equation}
W(\om,\bk,\bx,0) = \frac{\sigma_s^2 T_s^2 \ell_s^{2d}}{4 k_o}
\big|\hat S(T_s \om,\ell_s \bk)\big|^2 \delta(\bx) .
\label{eq:paraxini}
\end{equation}

Note that the right-hand side of equation \eqref{eq:TRresultP} is a
convolution, so we can write the Wigner transform explicitly using
Fourier transforms, as explained in Appendix~\ref{ap:SolPar}. The
result is
\begin{align}
W(\om,\bk,&\bx,z) = \frac{\sigma_s^2 T_s \ell_s^d}{4 k_o} \int_{\RR}
\frac{{\rm d} \Omega}{2 \pi} \int_{\RR^d} \frac{{\rm d} \bK}{(2
  \pi)^d} \, |\hat S(\Omega, \bK)|^2 \int_{\RR}{\rm d} t \int_{\RR^d}
     {\rm d} \by \nonumber \\&\times \int_{\RR^d} \frac{{\rm d}
       \bq}{(2 \pi)^d}\, \exp \left\{i\Big(\om-\frac{\Omega}{T_s}
     \Big)t -i \by \cdot \Big(\bk - \frac{\bK}{\ell_s}\Big)+i \bq
     \cdot \Big(\bx -\frac{\bK}{k_o} \frac{z}{\ell_s} \Big)\right.
     \nonumber \\ & \left. \hspace{0.3in}+ \frac{\sigma_c^2 \ell
       k_o^2}{4} \int_0^z {\rm d} z' \, \Big[
       \cR\Big(\frac{t}{T},\frac{\by - \frac{\bq}{k_o}(z- z') -
         \bv_o t}{\ell} \Big) - \cR(0,{\bf 0})\Big]
     \right\},\label{eq:Sol0}
\end{align}
and we use it next in the inverse problem of estimating the source
location and the mean flow velocity $\vv_o$.
\section{Application to imaging}
\label{sect:Inverse}
In this section we use the transport theory in the paraxial regime,
stated in section \ref{sect:resPar}, to localize a stationary in space
time-harmonic source in a moving random medium with smooth and
isotropic random fluctuations, from measurements at a stationary array
of receivers.  
The case of a time-harmonic source is interesting because it shows the
beneficial effect of the motion of the random medium for imaging. In
the absence of this motion, the wave received at the array is
time-harmonic, it oscillates at the frequency $\om_o$, and it is
impossible to determine from it the range of the source.  The random
motion of the medium causes broadening of the frequency support of the
wave field, which makes the range estimation possible.

We consider a strongly scattering regime, where the wave received at
the array is incoherent.  This means explicitly that the range $L$ is
much larger than the scattering mean free path $\mathscr{S}_{\rm par}$
or, equivalently, from \eqref{eq:TrP3},
\begin{equation}
\frac{\sigma_c^2 \ell k_o^2 L}{4} \cR(0,{\bf 0}) \gg 1.
\label{eq:StrScat}
\end{equation}
We also suppose that 
\begin{equation}
\frac{\eta}{\eta_s} = \frac{T}{T_s} \ll 1,
\label{eq:stability}
\end{equation}
to ensure that the imaging functions are statistically stable with
respect to the realizations of the random medium.
We begin in section \ref{sect:Inv1} with the approximation of the
Wigner transform \eqref{eq:Sol0} for a time-harmonic source, in the
strongly scattering regime. This Wigner transform quantifies the
time-space coherence properties of the wave, as described in section
\ref{sect:TXcoh}.  Then, we explain in section \ref{sect:Inv2} how we
can estimate the Wigner transform from the measurements at the
array. The source localization problem is discussed in section
\ref{sect:Inv3} and the estimation of the mean medium velocity is
discussed in section \ref{sect:Inv4}.

\subsection{Wigner transform for time-harmonic source and strong scattering} 
\label{sect:Inv1}
To derive the Wigner transform for a time-harmonic source, we take the
limit $T_s \to \infty$ in \eqref{eq:Sol0}, after rescaling the source
amplitude as
\begin{equation}
\sigma_s = \sigma/\sqrt{T_s}, \quad \sigma = O(1).
\label{eq:I1}
\end{equation}
We assume for convenience\footnote{The results extend qualitatively to
  other profiles but the formulas are no longer explicit.}  that the
source has a Gaussian profile,
\begin{equation}
\int_{\RR} {\rm d} \Omega \, \big|\hat S(\Omega,\bK)\big|^2 = (2
\pi)^{d} e^{-|\bK|^2},
\label{eq:I2}
\end{equation}
so we can calculate explicitly the integral over $\bK$ in
\eqref{eq:Sol0}. We obtain after the change of variables $ \by = \bxi
+ ({\bq}/{k_o}) z,$ that
\begin{align}
W(\om,\bk,\bx,z) =& \frac{\sigma^2\ell_s^d \pi^{d/2}}{4 k_o(2 \pi)^d}
\int_{\RR} {\rm d} t \, \int_{\RR^d} {\rm d} \bxi \int_{\RR^d} {\rm d}
\bq \, \exp \left\{ i \om t - \frac{|\bxi|^2}{4 \ell_s^2} - i \bxi
\cdot \bk + i \bq \cdot \Big(\bx -\frac{\bk}{k_o} z
\Big)\right. \nonumber \\ &\hspace{0.3in}\left. + \frac{\sigma_c^2
  \ell k_o^2}{4} \int_0^z {\rm d} z' \, \Big[
  \cR\Big(\frac{t}{T},\frac{\bxi+\frac{\bq}{k_o}z' - \bv_o t}{\ell}
  \Big) - \cR(0,{\bf 0})\Big] \right\}.\label{eq:I3}
\end{align}

Note that the last term in the exponent in \eqref{eq:I3} is negative,
because $\cR$ is maximal at the origin. Moreover, the relation
\eqref{eq:StrScat} that defines the strongly scattering regime implies
that the integrand in \eqref{eq:I3} is negligible for $t/T \ge 1$ and
$|\bxi+\bq/k_o z' - \bv_o t| /\ell \ge 1$.  Thus, we can restrict the
integral in \eqref{eq:I3} to the set
\[
\Big\{(t,\bxi,\bq) \in \RR^{2d+1}: ~ ~ |t| \ll T,~
~ \big|\bxi +\frac{\bq}{k_o}z' - \bv_o t \big |  \ll \ell \Big\},
\]
and approximate
\begin{equation}
\cR(\tau,\br) \approx \cR(0,{\bf 0}) - \frac{\alpha_o}{2} \tau^2
-\frac{\vartheta_o}{2} |\br|^2,
\label{eq:Taylor}
\end{equation}
with $\alpha_o, \vartheta_o > 0$. Here we used that the Hessian of
$\cR$ evaluated at the origin is negative definite and because the
medium is statistically isotropic, it is also diagonal, with the
entries $-\alpha_o$ and $-\vartheta_o$.  We obtain that
\begin{equation*}
\frac{\sigma_c^2 \ell k_o^2}{4}
\left[\cR\Big(\frac{t}{T},\frac{\bxi+\frac{\bq}{k_o}z' - \bv_o
    t}{\ell} \Big) - \cR(0,{\bf 0})\right] \approx -\frac{\alpha}{2}
\Big(\frac{t}{T}\Big)^2 - \frac{\vartheta}{2}
\bigg(\frac{|\bxi+\frac{\bq}{k_o}z' - \bv_o t|}{\ell} \bigg)^2
\end{equation*}
with the positive parameters 
\begin{equation}
\label{eq:pn}
 \alpha = \alpha_o \frac {\sigma_c^2 \ell k_o^2}{4} , \quad \quad
 \vartheta = \vartheta_o \frac{\sigma_c^2 \ell k_o^2}{4}.
\end{equation}
Substituting in \eqref{eq:I3} and integrating in $z'$ we obtain
\begin{align}
W(\om,\bk,\bx,z) &\approx \frac{\sigma^2 \ell_s^d \pi^{d/2}}{4 k_o(2
  \pi)^d} \int_{\RR} {\rm d} t \, \int_{\RR^d} {\rm d} \bxi
\int_{\RR^d} {\rm d} \bq \, \exp \left\{ i \om t - \frac{\alpha z}{2}
\Big(\frac{t}{T}\Big)^2 -\frac{\vartheta z}{2 \ell^2} \big|\bxi- \bv_o
t\big|^2 \right. \nonumber \\ &\hspace{-0.2in}\left. -
\frac{|\bxi|^2}{4 \ell_s^2}- i \bxi \cdot \bk - \frac{\vartheta z^2}{2
  \ell^2}(\bxi - \bv_o t)\cdot \frac{\bq}{k_o} - \frac{\vartheta
  z^3}{6 \ell^2} \Big| \frac{\bq}{k_o}\Big|^2+i \bq \cdot \Big(\bx
-\frac{\bk}{k_o} z\Big) \right\}.\label{eq:IM3}
\end{align}
The imaging results are based on this expression.  Before we present
them, we study the coherence properties of the transmitted wave and
define the coherence parameters which affect the performance of the
imaging techniques.

\subsection{Time-space coherence}
\label{sect:TXcoh} 
Let us define the time-space coherence function 
\begin{equation}
\label{eq:Coh1}
 {C}(\Delta t, \Delta \bx,\bx,z) = \frac{\lambda_o}{2\pi (c_o
     \rho_o)^2}\int_\RR {\rm d}t \, p(t+\Delta t,\bx+\Delta\bx/2,z)
   \overline{ p(t,\bx-\Delta\bx/2,z) } e^{i \omega_o \Delta t},
\end{equation}
and obtain from \eqref{eq:B3} that in the paraxial regime
\begin{align}
C(\Delta t, \Delta \bx,\bx,z) \approx& \int_{\RR } \frac{{\rm d}
  \omega}{2\pi} \int_{\RR^{d }} \frac{ {\rm d}\bk}{(2\pi)^d}
\int_{\RR^{d }} \frac{{\rm d}\bq}{(2\pi)^d} \, a(\omega,\bk+\bq/2,z)
\overline{a(\omega,\bk-\bq/2,z)} \nonumber \\ &\times
\exp\Big\{i \bq \cdot [z \nabla \beta(\bk) + \bx] +i
\Delta \bx\cdot \bk -i \omega \Delta t\Big\}. \label{eq:Coh2}
\end{align}
Moreover, in view of \eqref{eq:T6:0} and the fact that we average in
time so that the statistical fluctuations of $C$ are small (see Remark
\ref{rem:stability}) we have
\begin{align}
C(\Delta t, \Delta \bx,\bx,z) &\approx \EE\Big[{C}(\Delta t, \Delta
  \bx,\bx,z)\Big] \nonumber \\&\approx \int_{\RR } \frac{{\rm d} \omega}{2\pi}
\int_{\RR^{d }} \frac{ {\rm d}\bk}{(2\pi)^d} W(\omega, \bk,\bx,z)
e^{-i \omega \Delta t+i \Delta \bx\cdot \bk } . 
\label{eq:Coh3}
\end{align}
This shows formally that we can characterize the Wigner transform as the
Fourier transform of the coherence function
\begin{equation}
  W(\omega, \bk,\bx,z) \approx \int_{\RR } {\rm d}\Delta t
  \int_{\RR^{d }} {\rm d}\Delta\bx \, {C}(\Delta t, \Delta
  \bx,\bx,z) e^{i \omega \Delta t-i \Delta \bx\cdot \bk } .
\label{eq:Coh4}
\end{equation}

Using the expression \eqref{eq:I3} of the Wigner transform in
\eqref{eq:Coh3} we find after evaluating the integrals that
\begin{align}
C(\Delta t, \Delta \bx,\bx, z) &\approx
 \frac{\sigma^2 \ell_s^d }{2^{2+d/2} k_o
    R_z^d } \, \exp \big[i \varphi( \Delta t, \Delta \bx,\bx, z)]
\nonumber \\ & \times \exp\bigg[ - \frac{\Delta t^2}{2{\cal T}_z^2} -
  \frac{|\bx|^2}{2 R_z^2} - \frac{|\Delta \bx|^2}{2 {\mathcal
      D}_{1z}^2} - \frac{|H_z \Delta \bx - \bv_o \Delta t |^2}{2
    {\cal D}_{2z}^2} \bigg] , \label{eq:Cbar}
 \end{align}
 with phase 
\begin{equation}
 \varphi( \Delta t, \Delta \bx,\bx, z) = \frac{k_o \bx\cdot \Big[
     \Big(1 + \vartheta z\Big(\frac{\ell_s}{\ell}\Big)^2\Big) \Delta
     \bx - \vartheta z\Big(\frac{\ell_s}{\ell}\Big)^2 \bv_o \Delta t
     \Big]}{z \Big[ 1 + \frac{2}{3}\vartheta
     z\Big(\frac{\ell_s}{\ell}\Big)^2\Big] } ,
\end{equation}
and coefficients
\begin{align}
{\cal T}_z =& \frac{T}{\sqrt{\alpha z}},  \quad \quad  R_z =
\frac{z}{\sqrt{2}\ell_s k_o} \left( 1 + \frac{2\ell_s^2}{3{\cal D}_z^2}
\right)^{1/2} , 
 \quad  \quad
{\cal D}_z =   \frac{\ell}{\sqrt{\vartheta z}}, 
\label{eq:Adef}
\\
{\cal D}_{1z} =& 2 {\cal D}_z \left[3 \Big(1+
  \frac{\ell_s^2}{6{\cal D}_z^2}\Big)\right]^{1/2} 
  ,\label{eq:Coeff1} \quad 
{\cal D}_{2z} = {\cal D}_z\left(\frac{1+ \frac{2\ell_s^2}{3{\cal D}_z^2}}{1+
  \frac{\ell_s^2}{6{\cal D}_z^2}}\right)^{1/2} , \quad  
   H_z = 1 - \frac{1}{2
  \Big(1+ \frac{\ell_s^2}{6 {\cal D}_z^2}\Big)} . 
  \end{align}

The decay of the coherence function in $\Delta x$ models the spatial
decorrelation of the wave on the length scale corresponding to the
characteristic speckle size. This is quantified by the length scales
${\cal D}_{1z}$ and ${\cal D}_{2z}$, which are of the order of ${\cal
  D}_z$. We call ${\cal D}_z$ the decoherence length and obtain from
\eqref{eq:F12} and \eqref{eq:pn} that it is of the order of the 
typical size $\ell$ of the random fluctuations of the medium,
\begin{equation}
\label{eq:defDz}
{\cal D}_z = \frac{\ell}{\sqrt{\vartheta z}} =
\frac{\ell}{\pi\sqrt{\vartheta_o}} \sqrt{\frac{L}{z}} = O(\ell).
\end{equation} 

The decay of the coherence function in $\Delta t$ models the temporal
decorrelation of the wave, on the time scale 
\begin{equation}
\mathcal{T}_z = \frac{T}{\sqrt{\alpha z}} = \frac{T}{\pi
  \sqrt{\alpha_o}}\sqrt{\frac{L}{z}} = O(T),
\label{eq:defTz}
\end{equation}
where we used definitions \eqref{eq:F12} and \eqref{eq:pn}. We call
${\cal T}_z$ the decoherence time and note that it is of the 
order of the life span $T$ of the random fluctuations of the
medium.

The decay of the coherence function in $|\bx|$ means that the waves
propagate in a beam with radius $R_z$, which evolves in $z$ as
described in \eqref{eq:Adef} and satisfies
\begin{equation}
R_z \approx \left\{ \begin{array}{lll}
\displaystyle  \frac{z}{\sqrt{2} \ell_s k_o} &
  \hbox{for} & \ell_s \ll {\cal D}_z,  \\
\displaystyle   \sqrt{\frac{\vartheta}{3}}
  \frac{z^{3/2}}{k_o \ell} & \hbox{for} & \ell_s \gg {\cal D}_z.
 \end{array}     
     \right.
\label{eq:defRz}
\end{equation}
This shows that the transition from diffraction based beam spreading
to scattering based beam spreading happens around the critical
propagation distance
\begin{equation}
     z^* = \frac{1}{\vartheta} \Big(\frac{\ell}{\ell_s} \Big)^2 = L
     \Big( \frac{\gamma_s}{\gamma} \Big)^2 \frac{1}{\pi^2 \vartheta_o
     }. \label{eq:defzstar}
\end{equation}
This expression is derived from equation $\ell_s = {\cal D}_{z^\star}$ and definitions
\eqref{eq:F12} and \eqref{eq:pn}, and it shows that $z^\star/L$ is finite in our
 regime\footnote{Recall from
  section \ref{sect:resPar} that the paraxial regime is obtained in
  the limit $\gamma \to 0$ so that $\gamma/\gamma_s = \ell_s/\ell$
  remains finite.  Here we allow the ratio $\ell_s/\ell$ to be large
  or small, but independent of $\gamma$ which tends to zero.}.

Note that when $z \gg z^*$ i.e., ${\cal D}_z \ll \ell_s$, the
coefficients \eqref{eq:Coeff1} become
\begin{align}
\label{eq:Coeff1Z}
 {\cal D}_{1z} \approx & \sqrt{2} \ell_s, \quad {\cal D}_{2z} \approx
 2 {\cal D}_z , \quad H_z \approx 1,
\end{align}
and the coherence function satisfies
\begin{equation}
 \frac{|C(\Delta t, \Delta \bx,\bx, z) |}{|C(0,{\bf 0},{\bf 0},z)|}
 \approx \exp\bigg( - \frac{\Delta t^2}{2{\cal T}_z^2} -
 \frac{|\bx|^2}{2 R_z^2 } - \frac{|\Delta \bx|^2}{4
   \ell_s^2} - \frac{| \Delta \bx - \bv_o \Delta t |^2}{8{\cal D}_z^2 }
 \bigg).
 \end{equation}
 Thus, the spatial spreading and decorrelation of the wave field for
 $z \gg z^*$ are governed by the parameters $R_z$, $\ell_s$ and ${\cal D}_z$,
 with $R_z$ given by the second case in \eqref{eq:defRz} and ${\cal D}_z$
 given in \eqref{eq:defDz}.  These parameters scale with the
 propagation distance $z<L$ as $R_z \sim z^{3/2}$ and ${\cal D}_z \sim
 z^{-1/2}$. The temporal decorrelation is on the scale ${\cal T}_z \sim z^{-1/2}$.
  
 \subsection{Estimation of the Wigner transform}
\label{sect:Inv2}
Suppose that we have a receiver array centered at $(\bx_o,z)$, with
aperture in the cross-range plane modeled by the appodization function
\begin{equation}
\mathscr{A}(\bx) = \exp\left(-\frac{|\bx-\bx_o|^2}{2
  (\varkappa/k_o)^2}\right).
\label{eq:EW1}
\end{equation}
The linear size of the array is modeled by the standard deviation $
\varkappa /k_o $, with dimensionless $\varkappa >0$
defining the diameter of the array expressed in units of $\la_o$.

Recalling the wave decomposition \eqref{eq:B3} and that $\beta(\bk)
\sim k_o$ in the paraxial regime, we define the estimated mode
amplitudes by
\begin{align}
a_{\rm est}(\om,\bk,z) &= \frac{k_o e^{-i \beta(\bk) z}}{i \om_o
  {\rho_o} } \int_\RR dt \int_{\RR^d} {\rm d} \bx \, \mathscr{A}(\bx)
p(t,\bx,z) e^{i(\om+\om_o)t} e^{-i \bk \cdot \bx} \nonumber \\ &=
\left( \frac{\varkappa^2}{2 \pi k_o^2}\right)^{d/2} \int_{\RR^d} {\rm
  d} \tilde \bk \, a(\om,\bk+\tilde \bk,z) e^{ i [\beta(\bk+\tilde
    \bk) - \beta(\bk)] z + i \tilde \bk \cdot \bx_o -
  \frac{\varkappa^2 |\tilde \bk|^2}{2 k_o^2}} .
\label{eq:EW2}
\end{align} 
With these amplitudes we calculate the estimated Wigner transform 
\begin{align}
W_{\rm est}(\om,\bk,\bx,z) &= \int_{\RR^d} \frac{{\rm d} \bq}{(2
  \pi)^d} e^{i \bq \cdot (\nabla \beta(\bk) z + \bx)} a_{\rm
  est}\Big(\om,\bk + \frac{\bq}{2},z\Big) \overline{ a_{\rm
    est}\Big(\om,\bk - \frac{\bq}{2},z\Big)} 
\end{align}
and obtain after  carrying out the
integrals and using the approximation
\begin{equation*}
 \Big[ \beta\Big(\bk + \frac{\bq}{2}\Big) - \beta\Big(\bk-
   \frac{\bq}{2}\Big)\Big]z \approx \bq \cdot \nabla \beta(\bk) z
 ,
\end{equation*}
that 
\begin{equation}
W_{\rm est}(\om,\bk,\bx,z) \approx
\Big(\frac{\varkappa^2}{\pi k_o^2}\Big)^{d/2}\hspace{-0.03in}
e^{-\frac{k_o^2|\bx-\bx_o|^2}{\varkappa^2}}
\int_{\RR^d}\hspace{-0.03in} {\rm d} \bK \, e^{-\frac{\varkappa^2
    |\bK|^2}{k_o^2}} W(\om,\bk + \bK,\bx,z) .
\label{eq:EW4}
\end{equation}  
We can now use the expression \eqref{eq:IM3} in this equation, to
obtain an explicit approximation for $W_{\rm est}$.  Equivalently, we
can substitute \eqref{eq:Coh4} in \eqref{eq:EW4} and obtain after
integrating in $\bK$ that
\begin{align}
W_{\rm est}(\om,\bk,\bx,z) \approx  e^{
  -\frac{k_o^2|\bx-\bx_o|^2}{\varkappa^2}} \hspace{-0.03in}\int_{\RR }
{\rm d}\Delta t \int_{\RR^{d }} \hspace{-0.03in} {\rm d}\Delta\bx \,
C(\Delta t, \Delta \bx,\bx,z) e^{i \omega \Delta t -i \Delta \bx\cdot
  \bk -\frac{k_o^2|\Delta \bx |^2}{4 \varkappa^2}} ,
\label{eq:EW5}
\end{align}  
with $C$ given in \eqref{eq:Cbar}. 

\begin{remark}
\label{rem:stability}
{Note from \eqref{eq:Coh1} that the time integration that defines the
coherence function is over a time interval determined by the pulse
duration $T_s$, which is larger than the coherence time $T$ of
the medium by assumption \eqref{eq:stability}.  
If we interpret the wave as a train of $T_s/T$ pulses of total duration $T$, each individual pulse
travels through uncorrelated layers of medium because 
the correlation radius of the medium $\ell$
is much smaller than $c_oT$. This follows from the fact that $\ell/(c_oT) = \eps/(\eta\gamma)$ and  $\eps \ll \gamma\eta$.
Thus, $C(\Delta t,
\Delta \bx,\bx,z)$ is the superposition of approximately $T_s/T$
uncorrelated components and its statistical fluctuations are small by
the law of large numbers. } Moreover, we conclude from \eqref{eq:EW5}
that the estimated Wigner transform is approximately equal to its
expectation, up to fluctuations of relative standard deviation that is
smaller than $\sqrt{T/T_s}$.
\end{remark}

\subsection{Source localization}
\label{sect:Inv3}
We now show how we can use the estimated Wigner transform to localize
the source. Recall that we use the system of coordinates with origin
at the center of the source. Thus, the location $(\bx_o,z)$ of the
center of the array relative to the source is unknown and the goal of
imaging is to estimate it.  We begin in section \ref{sect:imagCR} with
the estimation of the direction of arrival of the waves at the array,
and then describe the localization in range in section
\ref{sect:imagR}. These two estimates determine the source location in
the cross-range plane, as well.

\subsubsection{Direction of arrival estimation}
\label{sect:imagCR}
We can estimate the direction of arrival of the waves from the peak
(maximum) in $\bk$ of the imaging function
\begin{align}
\mathcal{O}_{\rm DoA}(\bk,z) &= \int_{\RR} \frac{{\rm d} \om }{2\pi}
W_{\rm est}(\om,\bk,\bx_o,z), 
\label{eq:DoA1}
\end{align}
determined by the estimated Wigner transform at the center of the
array of receivers. If the medium were homogeneous, the maximum of
$\bk \mapsto \mathcal{O}_{\rm DoA}(\bk,z)$ would be at the cross-range
wave vector $\bk^* = k_o \frac{\bx_o}{z},$ and the width of the peak
(the resolution) would be $1/(\sqrt{2} \varkappa)$.  However,
cumulative scattering in the random medium gives a different result,
as we now explain:

Substituting \eqref{eq:EW5} in \eqref{eq:DoA1} and using the
expression \eqref{eq:Cbar}, we obtain after evaluating the integrals
that
\begin{align} 
\frac{ \mathcal{O}_{\rm DoA}(\bk,z) } { \max_{\bk'} \mathcal{O}_{\rm
    DoA}(\bk',z) } = \exp \left\{- \frac{1}{2 \vartheta_{\rm
    DoA}^2(z)} \left| \frac{\bk - \bk(z)}{k_o}
\right|^2\right\},
 \label{eq:DoA5}
\end{align}
  with
\begin{equation}
\vartheta_{\rm DoA}(z) = \left\{ \frac{1}{3 {\cal D}_z^2 k_o^2} \left(\frac{1
    + \frac{\ell_s^2}{2 {\cal D}_z^2}}{ 1 + \frac{2 \ell_s^2}{3 {\cal D}_z^2}}
    \right) + \frac{1}{2 \varkappa^2} \right\}^{1/2},
    \quad \quad
    \bk(z) = k_o \frac{\bx_o}{z} \left(\frac{1 + \frac{\ell_s^2}{{\cal D}_z^2}}{
  1 + \frac{2\ell_s^2}{3{\cal D}_z^2} }\right).
\label{eq:DoA6}
\end{equation}
Therefore, the maximum of $\bk \mapsto \mathcal{O}_{\rm DoA}(\bk,z)$ 
is at the cross-range wave vector $\bk(z)$
and the width of the peak (the resolution) is determined by
$\vartheta_{\rm DoA}(z)$. This resolution improves for larger array aperture
(i.e., $\varkappa$) and deteriorates as $z $ increases.  Depending on
the magnitude of $z$ relative to the critical range $z^*$ defined in 
\eqref{eq:defzstar}, we distinguish three cases:

1. In the case $z \ll z^*$ i.e., $\ell_s \ll {\cal D}_z$,  the intensity
travels along the deterministic characteristic, meaning that
$\mathcal{O}_{\rm DoA}(\bk,z)$ peaks at
\begin{equation}
\bk(z) \approx k_o \frac{\bx_o}{z}.
\end{equation}
However, the resolution is worse than in the homogeneous medium, 
\begin{equation}
\vartheta_{\rm DoA}(z) \approx \left\{ \frac{1 }{3 {\cal D}_z^2 k_o^2} +
\frac{1}{2\varkappa^2} \right\}^{1/2},
\label{eq:DoA6a}
\end{equation} 
with ${\cal D}_z$ defined in \eqref{eq:defDz}.

2. In the case $z \gg z^*$, i.e., $\ell_s \gg {\cal D}_z$, the peak of
$\mathcal{O}_{\rm DoA}(\bk,z) $ is at the cross-range wave vector
\begin{equation}
\bk(z) \approx \frac{3}{2} k_o \frac{\bx_o}{z},
\end{equation}
and the resolution is  
\begin{equation}
\vartheta_{\rm DoA}(z) \approx \left\{ \frac{1 }{4 {\cal D}_z^2 k_o^2} +
\frac{1}{2\varkappa^2} \right\}^{1/2}.
\label{eq:DoA6b}
\end{equation} 
Here the peak corresponds to a straight line characteristic, but with
a different slope than in the homogeneous medium. The resolution is
also worse than in the homogeneous medium.

3. In the case $z=O(z^*)$, the characteristic can no longer be
approximated by a straight line, as seen from
\eqref{eq:DoA6}. Nevertheless, we can still estimate the source
position from the observed peak $\bk(z)$, provided that we have an
estimate of the range~$z$.  The resolution of the estimate of $\bk(z)$
is $\vartheta_{\rm DoA}(z) $ given by (\ref{eq:DoA6}) that is bounded
from below by (\ref{eq:DoA6b}) and from above by (\ref{eq:DoA6a}).

\vspace{0.05in}
\begin{remark}
\label{rem:DOA}
Note that (\ref{eq:DoA6}) is a decreasing function of the array
diameter $\varkappa /k_o$, as long as this satisfies $\varkappa /k_o
\le \sqrt{2} {\cal D}_z$. Thus, increasing the aperture size beyond the
critical value $\sqrt{2} {\cal D}_z$ does not bring any resolution
improvement.
\end{remark}
\subsubsection{Range estimation}
\label{sect:imagR}
The results of the previous section show that the direction of arrival
estimation is coupled with the estimation of the range $z$ in general,
with the exception of the two extreme cases 1.~and 2.~outlined above.

To estimate the range $z$, we use the imaging function
\begin{align}
\mathcal{O}_{\rm range}(t,z) &= \int_{\RR} \frac{{\rm d} \om }{2\pi}
e^{-i \om t} \int_{\RR^d} \frac{{\rm d} \bk}{(2\pi)^d} W_{\rm
  est}(\om,\bk,\bx_o,z)  \approx {C}(t,{\bf 0},\bx_o,z),
\label{eq:Ra1}
\end{align}   
derived from \eqref{eq:EW5}. Substituting the expression
\eqref{eq:Cbar} of the coherence function in this equation we obtain
\begin{align}
\frac{ |\mathcal{O}_{\rm range}(t,z)| } { \max_{t'} |\mathcal{O}_{\rm
    range}(t',z)| } &= \exp \left\{- \frac{t^2}{2 \vartheta_{\rm
    range}^2(z)} \right\},
\label{eq:Ra2}
\end{align}  
with
\begin{align}
\vartheta_{\rm range}(z) = {\cal T}_z \left\{ 1 + \frac{|\bv_o|^2
  {\cal T}_z^2}{{\cal D}_z^2} \left( \frac{1 + \frac{\ell_s^2}{6 {\cal
      D}_z^2}}{1 + \frac{2 \ell_s^2}{3 {\cal
      D}_z^2}}\right)\right\}^{-1/2}.
\label{eq:Ra3}
\end{align}
As a function of $t$, this peaks at $t = 0$ and its absolute value
decays as a Gaussian, with standard deviation $\vartheta_{\rm
  range}(z)$.  If we know the statistics of the medium (the
decoherence time ${\cal T}_z$ and length ${\cal D}_z$) and the
magnitude of the cross-range velocity $|\bv_o|$, then we can determine
the range $z$ by estimating the rate of decay of $\mathcal{O}_{\rm
  range}(t,z)$. Note that the array dimameter $\varkappa /k_o$ plays
no role for the range estimation.

 
\begin{remark}
\label{rem:range2}
We can also estimate the mean velocity $\vv_o = (\bv_o,v_{oz})$ from
\eqref{eq:Ra2}, by considering different beam orientations in the case
that the sources and also the medium statistics (the decoherence
time ${\cal T}_z$ and length ${\cal D}_z$) are known.  That is to say,
with three known beams we can get the vector $\vv_o$, and then we can
use it to localize the unknown source using the direction of arrival
and range estimation described above. See also section \ref{sect:Inv4}
for a more detailed analysis of the velocity estimation.
\end{remark}
\begin{remark}
\label{rem:range3}
If the decoherence time ${\cal T}_z$ and length ${\cal D}_z$ are not
known, they can also be estimated using additional known
sources. Definitions \eqref{eq:defDz}--\eqref{eq:defTz} show that
${\cal D}_z z^{1/2}$ and ${\cal T}_z z^{1/2}$ are constant with
respect to $z$. Once estimated, these constants can be used in the
imaging of the unknown source.
\end{remark}
\subsection{Single beam lateral velocity estimation}\label{sect:Inv4}
We observe from \eqref{eq:DoA5} and \eqref{eq:Ra2}--\eqref{eq:Ra3}
that the source localization depends only on the Euclidian norm
$|\bv_o|$ of the cross-range component of the mean velocity of the
medium.  We show here that $\bv_o$ can be obtained with only one beam
and, when the receiver array is large and $z \gg z^*$ i.e., $\ell_s
\gg D_z$, the velocity estimate is independent of the medium
statistics and the source location.

The estimation of $\bv_o$ is based on the imaging function
\begin{align}
 \mathcal{O}_{\rm v}(\by,t,z) &= \int_{\RR}\frac{ {\rm d} \om}{2\pi} e^{-i
    \om t} \int_{\RR^d} \frac{{\rm d} \bk }{(2\pi)^d} e^{ i \bk \cdot
    \by} \int_{\RR^d} {\rm d} \bx \, W_{\rm est}(\om,\bk,\bx,z)
  \nonumber \\ &\approx \exp\Big( - \frac{k_o^2 |\by|^2}{4 \varkappa^2}
  \Big) \int_{\RR^d} {\rm d} \bx \, \exp\Big(
  -\frac{k_o^2|\bx-\bx_o|^2}{\varkappa^2}\Big) \, {C}(t,\by ,\bx, z) .
\end{align}
Substituting the expression \eqref{eq:Cbar} of the coherence function
and carrying out the integrals we obtain that
\begin{align}
 |\mathcal{O}_{\rm v}(\by,t,z)|  \approx &  \frac{\sigma^2 \pi^{d/2 } (\varkappa\ell_s)^d  }{2^{2+d} k_o^{d+1}   {\cal A}^{d}_z  }   \exp\left\{ - \frac{t^2}{2 {\cal
    T}_z^2} -\frac{|\by- s_z t \bv_o|^2}{2 \mathfrak{m}^2_z{\cal A}^2_z} -
\frac{|t \bv_o|^2}{ \mathfrak{n}^2_z {\cal A}_z^2} - \frac{| \bx_o|^2}{4
  {\cal A}^2_z} \right\}, 
\label{eq:Ovy}
\end{align}
with the effective apperture 
\begin{align}
 {\cal A}^2_z = \frac{1}{4} \Big[\Big(\frac{\varkappa}{k_o}\Big)^2 +
   \Big(\frac{z}{k_o \ell_s}\Big)^2 \Big(1 + \frac{2 \ell_s^2}{3
     D_z^2}\Big)\Big],
\label{eq:defAz}
\end{align} 
and dimensionless parameters
\begin{align*}
 {\mathfrak m}_z^2 =& \frac{8}{1 +
   \frac{2}{3 {\cal D}_z^2} \big(\frac{z}{k_o \ell_s}\big)^2\Big(1 +
   \frac{\ell_s^2}{2 {\cal D}_z^2}\Big) + \big(\frac{\varkappa}{k_o
     \ell_s}\big)^2 \Big(1 + \frac{2 \ell_s^2}{{\cal D}_z^2}\Big) +
   \big(\frac{k_o}{\varkappa}\big)^2 \big(\frac{z}{k_o \ell_s}\big)^2
   \Big( 1 + \frac{2 \ell_s^2}{3 {\cal D}_z^2}\Big)}, \\ \mathfrak{n}_z^2 =
 &  \frac{\mathfrak{m}_z^2}{s_z(q_z - s_z/2)},\\
s_z =& \frac{\mathfrak{m}_z^2}{4{\cal D}_z^2} \Big[
  \Big(\frac{\varkappa}{k_o}\Big)^2 + \frac{1}{2} \Big(\frac{z}{k_o
    \ell_s}\Big)^2 \Big(1 + \frac{\ell_s^2}{3 {\cal D}_z^2} \Big)\Big],
\quad \quad q_z = \frac{\Big(\frac{\varkappa}{k_o}\Big)^2 + 
 \Big(\frac{z}{k_o
    \ell_s}\Big)^2 \Big(1 + \frac{\ell_s^2}{6 {\cal D}_z^2} \Big)}{
2\Big(\frac{\varkappa}{k_o}\Big)^2+ 
\Big(\frac{z}{k_o
    \ell_s}\Big)^2 \Big(1 + \frac{\ell_s^2}{3 {\cal D}_z^2} \Big)}.
 \end{align*}
These depend on the radii $\varkappa/k_o$ of the array and $\ell_s$ of
the source, the decoherence length ${\cal D}_z$ and the ratio
$z/(k_o \ell_s)$ that quantifies the cross-range resolution of
focusing of a wave using time delay beamforming at a source of radius
$\ell_s$.

To estimate $\bv_o$ we can proceed as follows: First, we estimate for
each time $t$ the position $\by_{\rm max}(t)$ that maximizes $
\by\mapsto \mathcal{O}_{\rm v}(\by;t,z)$.  Second, we note from
(\ref{eq:Ovy}) that $\by_{\rm max}(t)$ should be a linear function in
$t$, of the form $\by_{\rm max}(t) = s_z \bv_o t$.  Therefore, we can
estimate $ s_z \bv_o$ with a weighted linear least squares regression
of $\by_{\rm max}(t) $ with respect to $t$. 
In practice $s_z$ is likely unknown. However, in the case of a large
receiver array with radius satisfying
\begin{equation}
\frac{\varkappa}{k_o} \gg \max\left\{\frac{z}{k_o \ell_s},
\frac{z}{k_o {\cal D}_z} \right\},
\label{eq:LA1}
\end{equation}
and for $z \gg z^\star$, so that $\ell_s \gg {\cal D}_z$, we obtain
that $s_z \approx 1$. Thus, the least squares regression gives an
unbiased estimate of $\bv_o$.
 
In view of (\ref{eq:Ovy}), the least squares regression can be carried
out over a time interval with length of the order of $\min(
{\cal T}_z , \mathfrak{n}_z {\cal A}_z/|{\itbf v}_o|) $. Beyond this
critical time the function $ \mathcal{O}_{\rm v}$ vanishes.  Therefore,
as long as
$|{\itbf v}_o|<  {{\mathfrak n}_z {\cal A}_z}/{{\cal T}_z} $,
the velocity resolution is 
\begin{align}
 \label{eq:resv1}
 \mbox{res}_{\rm v}
  = \frac{{\mathfrak m}_z {\cal A}_z}{s_z {\cal T}_z } \approx
  \frac{{\cal D}_z}{{\cal T}_z},
\end{align}
where the approximation is for a large array and $\ell_s \gg {\cal
  D}_z$.  
If $|{\itbf v}_o|$ is larger than $ {{\mathfrak n}_z {\cal A}_z}/{{\cal T}_z} $, then the resolution is reduced to
\begin{align} 
\label{eq:resv2}
\mbox{res}_{\rm v} = 
\frac{{\mathfrak m}_z {\cal A}_z}{s_z {\mathfrak n}_z {\cal
    A}_z/|\bv_o| } \approx \frac{{\cal D}_z}{{\cal T}_z} \frac{|\bv_o|
  {\cal T}_z}{\mathfrak{n}_z {\cal A}_z}.
\end{align}

\section{Analysis of the wave field}
\label{sect:Anal}
To derive the results stated in section \ref{sect:ResForward}, we
begin in section \ref{sect:Anal1} with a slight reformulation, which
transforms equation \eqref{eq:P2} into a form that is more convenient
for the analysis. We scale the resulting equation   in section
\ref{sect:Anal2}, in the regime defined in section \ref{sect:scaling},
and then we change coordinates to a moving frame in section
\ref{sect:Anal3}. In this frame we write the wave as a superposition
of time-harmonic, plane waves with random amplitudes that model the
net scattering in the random medium, as described in section
\ref{sect:Anal4}.  We explain in section \ref{sect:Anal5} that the
backward going waves are negligible, and use the diffusion
approximation theory in section \ref{sect:Anal6} to analyze the
amplitudes of the forward going waves, in the limit $\ep \to 0$.  We
end in section \ref{sect:Anal7} with the paraxial limit.

\subsection{Transformation of the wave equation}
\label{sect:Anal1}
Let us define the new potential 
\begin{equation}
\psi(t,\vx) = \frac{\sqrt{\rho(t,\vx)}}{\sqrt{\rho_o}} \phi(t,\vx),
\label{eq:An1}
\end{equation}
and substitute it in \eqref{eq:P2} to obtain the wave equation
\begin{align}
D_t \left[ \frac{1}{c^2(t,\vx)} D_t \psi(t,\vx)\right] - \frac{D_t
  \psi(t,\vx)D_t \ln \rho(t,\vx)}{c^2(t,\vx)}
-\Delta_{\vx} \psi(t,\vx) +q(t,\vx) \psi(t,\vx) \nonumber \\ 
=\sigma_s \frac{\sqrt{\rho(t,\vx)}}{\sqrt{\rho_o}} e^{-i \om_o t} S
\Big(\frac{t}{T_s},\frac{\bx}{\ell_s}\Big) \delta(z),
\label{eq:An3}
\end{align}
for $t \in \RR$ and $\vx = (\bx,z) \in \RR^{d+1}$, where
$\Delta_{\vx}$ is the Laplacian operator and
\begin{align}
q(t,\vx) = \frac{\Delta_{\vx} \sqrt{\rho(t,\vx)}}{\sqrt{\rho(t,\vx)}}
- \frac{1}{c^2(t,\vx)} \left\{ \frac{D^2_t
  \sqrt{\rho(t,\vx)}}{\sqrt{\rho(t,\vx)}} - \frac{1}{2}\left[ D_t \ln
  \rho(t,\vx) \right]^2 \right\} \nonumber \\-\frac{1}{2}D_t
c^{-2}(t,\vx) D_t \ln{\rho(t,\vx)}.
\label{eq:An4}
\end{align}
The initial condition \eqref{eq:P3} becomes
\begin{equation}
\psi(t,\vx) \equiv 0, \quad t \ll -T_s.
\label{eq:An4p}
\end{equation}

\subsection{Scaled wave equation}
\label{sect:Anal2}
We use the scaling regime defined in section \ref{sect:scaling} and
denote with primes the dimensionless, order one variables
\begin{equation}
\vx = L \vx', \quad t = T_L t'.
\label{eq:An6}
\end{equation}
We also let
\begin{equation}
\vv_o = V \vv_o', \quad c_o = c_o c_o', \quad \om_o = \om_o
\frac{\om_o'}{2 \pi},
\label{eq:An7}
\end{equation}
where the constants $c_o' = 1$ and $\om_o' = 2 \pi$ are introduced so
that the scaled equation is easier to interpret.

In the scaled variables, and using the source amplitude \eqref{eq:F14}, 
the right-hand side in \eqref{eq:An3} becomes 
\begin{align}
\hspace{-0.1in}\sigma_s \frac{\sqrt{\rho(t,\vx)}}{\sqrt{\rho_o}} e^{-i \om_o t} S
\Big(\frac{t}{T_s},\frac{\bx}{\ell_s}\Big) \delta(z) = \frac{\big[1+
    O(\sqrt{\ep})\big]}{\ep \eta_s L^2} \Big(\frac{\gamma_s}{\ep}
\Big)^d e^{-i \frac{\om_o'}{\ep} t'} S
\Big(\frac{t'}{\eta_s},\frac{\gamma_s \bx'}{\ep} \Big) \delta (z').
\label{eq:An7p}
\end{align}
We also have from definitions \eqref{eq:F1}--\eqref{eq:F3} that the
random coefficients take the form
\begin{align}
\frac{\vv(t,\vx)}{V} &= \vv'(t',\vx') = \vv_o' + \sqrt{\ep \gamma}
\,\bar\sigma_{v} \,\vnu \Big(\frac{t'}{\eta}, \frac{\vx'-\ep \vv_o'
  t'}{\ep/\gamma}\Big),
\label{eq:An9} \\
\frac{\rho(t,\vx)}{\rho_o} &= \exp \left[ \sqrt{\ep \gamma}\, 
  \bar\sigma_{\rho} \, \nu_\rho \Big(\frac{t'}{\eta}, \frac{\vx'-\ep
    \vv_o't'}{\ep/\gamma} \Big) \right], \label{eq:An9p}
\\ \frac{c_o^2}{c^2(t,\vx)} &= \frac{1}{(c_o')^2} \left[ 1 + \sqrt{\ep
    \gamma}\, \bar{\sigma}_c \,\nu_c\Big(\frac{t'}{\eta}, \frac{\vx'-\ep
    \vv_o't'}{\ep/\gamma} \Big) \right], \label{eq:An9pp}
\end{align}
with scaled standard deviations $ \bar{\sigma}_c,\bar{\sigma}_v,\bar{\sigma}_\rho$  defined in \eqref{eq:F12}.

The solution $\psi$ of \eqref{eq:An3}
must have variations in $t'$ and $\vx'$ on the same scale as the
source term and the coefficients \eqref{eq:An9}--\eqref{eq:An9pp}, meaning that
$\partial_{t'} \psi \sim  {{1}/{\ep}}$, $|\nabla_{\vx'} \psi| \sim
{1}/{\ep}$.
From equations \eqref{eq:F8}--\eqref{eq:F13} we obtain that in the
scaled variables we have 
\begin{equation}
D_t = \frac{1}{T_L c_o'} D_{t'}^\ep, \quad \mbox{with}~D_{t'}^\ep =
\partial_{t'} + \ep \vv'(t',\vx') \cdot \nabla_{\vx'}.
\label{eq:An8}
\end{equation}
Equation \eqref{eq:An9p} gives
\begin{align*}
D_t \ln{\rho(t,\vx)} &= \frac{\sqrt{\ep
    \gamma} \bar{\sigma}_\rho}{ T_L } D_{t'}^\ep \nu_\rho
\Big(\frac{t'}{\eta}, \frac{\vx'-\ep \vv_o't'}{\ep/\gamma} \Big) =
\frac{O(\sqrt{\ep})}{T_L} ,
\end{align*}
and  
\begin{align*} 
\frac{D_t^2 \sqrt{\rho(t,\vx)}}{
  \sqrt{\rho(t,\vx)}} &   = \frac{\sqrt{\ep \gamma} \bar{\sigma}_\rho}{2
  T_L^2} (D_{t'}^\ep)^2 \nu_\rho \Big(\frac{t'}{\eta}, \frac{\vx'-\ep
  \vv_o't'}{\ep/\gamma} \Big) + \frac{\ep \gamma
  \bar{\sigma}_\rho^2}{4 T_L^2} \left[ D_{t'}^\ep \nu_\rho
  \Big(\frac{t'}{\eta}, \frac{\vx'-\ep \vv_o't'}{\ep/\gamma}
  \Big)\right]^2 \nonumber  \\&= \frac{O(\sqrt{\ep})}{T_L^2}. 
\end{align*}
From equation \eqref{eq:An9pp} we get
\begin{equation*}
D_t \left[ \frac{1}{c^2(t,\vx)} \right] = \frac{\sqrt{\ep
    \gamma} \bar{\sigma}_c}{(c_o c_o')^2 T_L} D_{t'}^\ep \nu_c \Big(\frac{t'}{\eta},
\frac{\vx'-\ep \vv_o't'}{\ep/\gamma} \Big) = \frac{O(\sqrt{\ep})}{c_o L},
\end{equation*}
and $q$ defined in \eqref{eq:An4} takes the form
\begin{equation}
q(t,\vx) = \frac{ \gamma^{5/2} \bar{\sigma}_\rho}{2 \ep^{3/2} L^2}\left[
Q^\ep\Big(\frac{t'}{\eta}, \frac{\vx'-\ep \vv_o't'}{\ep/\gamma}
\Big) + O(\ep^2)\right],
\label{eq:An15}
\end{equation}
with 
\begin{equation}
Q^\ep(\tau,\vr) = \Delta_{\vr} \nu_\rho(\tau,\vr) + \frac{\sqrt{\ep
    \gamma} \bar{\sigma}_\rho}{2}  \left|\nabla_{\vr} \nu_\rho(\tau,\vr)
\right|^2.
\label{eq:An16}
\end{equation}

Substituting in \eqref{eq:An3} and multiplying both sides by $\ep
L^2$, we obtain that  the potential denoted by
$\psi'(t',\vx')$ in the scaled variables satisfies 
\begin{align}
&\hspace{-0.15in}\ep \left\{ \frac{\left[1 + \sqrt{\ep \gamma}\,
      \bar{\sigma}_c \nu_c\big(\frac{t'}{\eta}, \frac{\vx'-\ep \vv_o't'}{\ep/\gamma}
      \big)\right]}{(c_o')^2} \partial_{t'}^2 + \frac{2 \ep}{(c_o')^2}
  \vv_o' \cdot \nabla_{\vx'} \partial_{t'} - \Delta_{\vx'}
  \right\}\psi'(t',\vx') \nonumber
  \\ & \hspace{-0.1in}+\frac{\bar{\sigma}_\rho \gamma^{5/2}}{2
    \sqrt{\ep}} Q^\ep\Big(\frac{t'}{\eta}, \frac{\vx'-\ep
    \vv_o't'}{\ep/\gamma} \Big) \psi'(t',\vx') \approx
  \frac{1}{\eta_s} \Big(\frac{\gamma_s}{\ep}\Big)^d e^{-i
    \frac{\om_o'}{\ep} t'} S
  \Big(\frac{t'}{\eta_s},\frac{\bx'}{\ep/\gamma_s} \Big) \delta (z'),
\label{eq:An18}
\end{align}
with initial condition obtained from \eqref{eq:P3} and \eqref{eq:An1},
\begin{equation}
\psi'(t',\vx') \equiv 0, \quad t' \ll - \eta_s.
\label{eq:An18p}
\end{equation}
The approximation in equation \eqref{eq:An18} is because we neglect
$O(\sqrt{\ep})$ terms that tend to zero in the limit $\ep \to 0.$ Note
in particular that the random perturbations $\vnu$ of the velocity of
the flow appear in these terms and are negligible in our regime.

All variables are assumed scaled 
in the remainder of the section and
we simplify notation by dropping the primes.
\subsection{Moving frame}
\label{sect:Anal3}
Let us introduce the notation $\vv_o = (\bv_o,v_{oz})$ for the scaled
mean velocity of the ambient flow, and change the range coordinate $z$
to
\begin{equation}
\zeta = z - \ep v_{oz} t.
\label{eq:An19}
\end{equation}
We denote the potential in this moving frame by
\begin{equation}
u^\ep(t,\bx,\zeta) = \psi(t,\bx,\zeta + \ep v_{oz} t), 
\label{eq:An20}
\end{equation}
and obtain from \eqref{eq:An18} that it satisfies the wave equation
\begin{align}
&\ep \left\{ \frac{\left[1 + \sqrt{\ep \gamma}\,
     \bar{\sigma}_c  \nu_c\Big(\frac{t}{\eta}, \frac{\bx-\ep \bv_ot}{\ep/\gamma},
      \frac{\gamma \zeta}{\ep} \Big)\right]}{c_o^2} \partial_{t}^2 +
   \frac{2 \ep}{c_o^2}   \bv_o \cdot \nabla_{\bx}   \partial_t   -
  \Delta_{\bx} - \partial_\zeta^2 \right\}u^\ep(t,\bx,\zeta) 
  \nonumber \\ & +\frac{\bar{\sigma}_\rho \gamma^{5/2}}{2
    \sqrt{\ep}}Q^\ep\Big(\frac{t}{\eta}, \frac{\bx-\ep
    \bv_ot}{\ep/\gamma},\frac{\gamma \zeta}{\ep} \Big)
  u^\ep(t,\bx,\zeta) \approx \Big(\frac{\gamma_s}{\ep}\Big)^d
  \frac{e^{-i \frac{\om_o}{\ep} t}}{\eta_s} S
  \Big(\frac{t}{\eta_s},\frac{\gamma_s \bx}{\ep} \Big) \delta (\zeta +
  \ep v_{oz} t),
\label{eq:An21}
\end{align}
where again we neglect the terms that become negligible in the limit
$\ep \to 0$. The gradient $\nabla_{\bx}$ and Laplacian $\Delta_{\bx}$
are in the cross-range variable $\bx \in \RR^d$.
\subsection{Wave decomposition}
\label{sect:Anal4}
The interaction of the waves with the random medium depends on the
frequency and direction of propagation, so we decompose
$u^\ep(t,\bx,\zeta)$ using the Fourier transform
\begin{equation}
\hat u^\ep(\om,\bk, \zeta) = \int_{\RR} {\rm d}t \int_{\RR^d} {\rm d} \bx \,
u^\ep(t,\bx,\zeta) e^{i \big(\frac{\om_o}{\ep} + \om\big)t - i
  \frac{\bk}{\ep} \cdot \bx},
\label{eq:An22}
\end{equation}
with inverse 
\begin{equation}
u^\ep(t,\bx, \zeta) = \int_{\RR} \frac{{\rm d} \om}{2 \pi}
\int_{\RR^d} \frac{{\rm d} \bk}{(2 \pi \ep)^d} \, \hat
u^\ep(\om,\bk,\zeta) e^{-i \big(\frac{\om_o}{\ep} + \om\big)t + i
  \frac{\bk}{\ep} \cdot \bx}.
\label{eq:An23}
\end{equation}
The Fourier transform of equation \eqref{eq:An21} is
\begin{align}
&\hspace{-0.15in}\left[ -\frac{\beta^2(\bk)}{\ep} - 2 k_o
    \Big(\frac{\om}{c_o} - \frac{\bv_o \cdot \bk}{c_o} \Big) - 
     \ep
    \partial_\zeta^2 \right] \hat u^\ep(\om,\bk,\zeta) 
- \frac{\eta
    \gamma^{1/2 -d}}{\sqrt{\ep}} \int_{\RR} \frac{{\rm d} \om'}{2\pi}
  \int_{\RR^d} \frac{{\rm d} \bk'}{(2\pi)^d} \nonumber
  \\ &\hspace{0.95in}\hat u^\ep(\om',\bk',\zeta) \left[k_o^2
   \bar{\sigma}_c \hat
    \nu_c\Big( \eta\big(\om-\om'-(\bk-\bk')\cdot \bv_o\big),
    \frac{\bk-\bk'}{\gamma}, \frac{\gamma \zeta}{\ep} \Big)
    \right. \nonumber \\ &\left. \hspace{1.2in} -
    \frac{\gamma^2\bar{\sigma}_\rho}{2} \hat Q^\ep \Big(
    \eta\big(\om-\om'-(\bk-\bk')\cdot \bv_o\big),
    \frac{\bk-\bk'}{\gamma}, \frac{\gamma \zeta}{\ep} \Big) \right]
  \nonumber \\ &\hspace{0.95in}\approx \frac{e^{-i \frac{\om \zeta}{\ep
        v_{oz}}}}{\ep \eta_s v_{oz}} \widecheck{S}\Big(
  -\frac{\zeta}{\ep \eta_s v_{oz}}, \frac{\bk}{\gamma_s}
  \Big), \label{eq:An25}
\end{align}
where 
\begin{equation}
\beta(\bk) = \sqrt{k_o^2 - |\bk|^2}, \quad k_o = \frac{\om_o}{c_o},
\label{eq:An24}
\end{equation}
and
\[
\widecheck{S}(\tau,\boldsymbol{\kappa}) = \int_{\RR^d} {\rm d}\br \,
S(\tau,\br) e^{-i \boldsymbol{\kappa} \cdot \br}.
\]

Note that the right-hand side in \eqref{eq:An25} is supported at
$|\bk| = O(\gamma_s)$, so by keeping $\gamma_s$ small, we ensure that
$\beta(\bk)$ remains real valued in our regime. Physically, this means
that $\hat u^\ep(\om,\bk,\zeta)$ is a propagating wave, not
evanescent.

Note also that if the mean velocity $\vv_o$ is orthogonal to the range
direction, the source term satisfies
\[
\lim_{v_{oz} \to 0} \frac{e^{-i \frac{\om \zeta}{\ep
        v_{oz}}}}{\ep \eta_s
    v_{oz}} \widecheck{S}\Big( -\frac{\zeta}{\ep \eta_s v_{oz}},
  \frac{\bk}{\gamma_s} \Big) \to \hat S\Big(\eta_s \om,
  \frac{\bk}{\gamma_s} \Big) \delta(\zeta),
\]
in the sense of distributions, where 
\[
\hat S(\omega,\boldsymbol{\kappa}) = \int_{\RR} {\rm d}\tau \int_{\RR^d}
     {\rm d}\br \, S(\tau,\br) e^{i \om \tau -i \boldsymbol{\kappa}
       \cdot \br}.
\]

We introduce
\begin{align}
 a^\ep(\om,\bk,\zeta)  & = \Big[ \frac{\sqrt{\beta(\bk)}}{2} 
\hat u^\ep(\om, \bk,\zeta)
+\frac{\ep}{2i  \sqrt{\beta(\bk)}} \partial_\zeta \hat u^\ep(\om, \bk,\zeta)
\Big] e^{- i \beta(\bk) \frac{\zeta}{\ep}},  \\
  a_-^\ep(\om,\bk,\zeta) 
  &= \Big[  \frac{\sqrt{\beta(\bk)}}{2} 
\hat u^\ep(\om, \bk,\zeta) -\frac{\ep}{2i  \sqrt{\beta(\bk)}} \partial_\zeta \hat u^\ep(\om, \bk,\zeta)
\Big] e^{ i \beta(\bk) \frac{\zeta}{\ep}} ,
\end{align}
so that we have the decomposition
\begin{align}
\hat u^\ep(\om, \bk,\zeta) = \frac{1}{\sqrt{\beta(\bk)}} \left[
  a^\ep(\om,\bk,\zeta) e^{i \beta(\bk) \frac{\zeta}{\ep}} +
  a_-^\ep(\om,\bk,\zeta) e^{-i \beta(\bk) \frac{\zeta}{\ep}}\right],
\label{eq:An26}
\end{align}
and the complex amplitudes $a^\ep$ and $a_-^\ep$ satisfy the relation
\begin{equation}
\partial_\zeta a^\ep(\om,\bk,\zeta) e^{i \beta(\bk) \frac{\zeta}{\ep}} +
\partial_\zeta a_-^\ep(\om,\bk,\zeta) e^{-i \beta(\bk) \frac{\zeta}{\ep}} = 0.
\label{eq:An27}
\end{equation}
This gives that 
\begin{align}
\partial_\zeta \hat u^\ep(\om, \bk,\zeta) = \frac{i
  \sqrt{\beta(\bk)}}{\ep} \left[ a^\ep(\om,\bk,\zeta) e^{i
    \beta(\bk) \frac{\zeta}{\ep}} - a_-^\ep(\om,\bk,\zeta) e^{-i \beta(\bk)
    \frac{\zeta}{\ep}}\right],
\label{eq:An28}
\end{align}
and moreover, that 
\begin{align}
\partial^2_\zeta \hat u^\ep(\om, \bk,\zeta) =
-\frac{\beta^2(\bk)}{\ep^2} \hat u^\ep(\om, \bk,\zeta) + \frac{2 i
  \sqrt{\beta(\bk)}}{\ep} \partial_\zeta a^\ep(\om,\bk,\zeta) e^{i
  \beta(\bk) \frac{\zeta}{\ep}}.
\label{eq:An29}
\end{align}

The  decomposition in \eqref{eq:An26} and   \eqref{eq:An23}    is a decomposition 
of  $u^\ep(t,\bx,\zeta)$ into  a superposition of plane waves  with wave
vectors 
\begin{equation} 
\vec{\bk}_{\pm} = (\bk,\pm \beta(\bk)),
\end{equation} where the plus sign denotes the
waves propagating in the positive range direction and the negative
sign denotes the waves propagating in the negative range direction. The amplitudes
$a^\ep$ and $a_-^\ep$  of these waves are random fields, which evolve in range
according to \eqref{eq:An27} and the equation
\begin{align}
&\partial_\zeta a^\ep(\om,\bk,\zeta) \approx \frac{i k_o (\om -\bv_o
    \cdot \bk)}{c_o\beta(\bk)} \left[ a^\ep(\om,\bk,\zeta) +
    a_-^\ep(\om,\bk,\zeta)e^{-2 i \beta(\bk) \frac{\zeta}{\ep}}
    \right] \nonumber \\ &\hspace{0.3in}+\frac{i \eta
    \gamma^{1/2-d}}{2 \sqrt{\ep}} \int_{\RR} \frac{{\rm d} \om'}{2\pi}
  \int \frac{{\rm d} \bk'}{(2 \pi)^d} \left[ k_o^2  \bar{\sigma}_c \hat \nu_c\Big(
    \eta\big(\om-\om'-(\bk-\bk')\cdot \bv_o\big),
    \frac{\bk-\bk'}{\gamma}, \frac{\gamma \zeta}{\ep} \Big)
    \right. \nonumber \\ &\hspace{1.55in}\left.
    -\frac{\gamma^2\bar{\sigma}_\rho}{2} \hat Q^\ep \Big(
    \eta\big(\om-\om'-(\bk-\bk')\cdot \bv_o\big),
    \frac{\bk-\bk'}{\gamma}, \frac{\gamma \zeta}{\ep} \Big)\right]
  \nonumber \\&\hspace{0.3in}\times \frac{1}{\sqrt{\beta(\bk)
      \beta(\bk')}}\left[a^\ep(\om',\bk',\zeta)e^{i
      \big[\beta(\bk')-\beta(\bk)\big] \frac{\zeta}{\ep}} +
    a_-^\ep(\om',\bk',\zeta)e^{-i \big[\beta(\bk') +
        \beta(\bk)\big]\frac{\zeta}{\ep}} \right] \nonumber
  \\ &\hspace{0.3in}+\frac{i}{2 \sqrt{\beta(\bk)} \ep \eta_s
    v_{oz}} \widecheck{S}\Big( -\frac{\zeta}{\ep \eta_s v_{oz}},
  \frac{\bk}{\gamma_s} \Big)e^{-i \frac{\om \zeta}{\ep v_{oz}}- i
    \beta(\bk) \frac{\zeta}{\ep}},
\label{eq:An30}
\end{align}
derived by substituting \eqref{eq:An26}--\eqref{eq:An29} into
\eqref{eq:An25}.

\subsection{Forward scattering approximation}
\label{sect:Anal5}
Equation \eqref{eq:An30} shows that the amplitudes $a^\ep$ are
coupled to each other and to $a^\ep_-$. In our scaling regime, where
the cone of directions of propagation has small opening angle
controlled by the parameter $\gamma_s$, and where the covariance
\eqref{eq:F5} of the fluctuations is smooth, the coupling between
$a^\ep$ and $a_{-}^\ep$ becomes negligible in the limit $\ep \to
0$. We refer to \cite[Section C.2]{BG2016} and \cite[Section
  5.2]{BG2016p} for a more detailed explanation of this fact.

Using the assumption that the random fluctuations are supported at
finite range (see section \ref{sect:form}), we require that the wave
be outgoing at $|\zeta| \to \infty$. This radiation condition and the
negligible coupling between $a^\ep$ and $a^\ep_-$ in the limit $\ep
\to 0$ imply that we can neglect the backward going waves, and we can
write
\begin{equation}
\hat u^\ep(\om,\bk,\zeta) \approx
\frac{a^\ep(\om,\bk,\zeta)}{\sqrt{\beta(\bk)}} e^{i \beta(\bk)
    \frac{\zeta}{\ep}}, \quad \zeta > O(\ep).
\label{eq:An31f}
\end{equation}

The starting value of $a^\ep(\om,\bk,\zeta)$ is determined by the
source term in \eqref{eq:An30}, which contributes only for $\zeta =
v_{oz}O(\ep)$. For such small $\zeta$, we can change variables $\zeta
= \ep \xi$ in \eqref{eq:An30} and obtain that
\begin{align*} 
\partial_\xi a^\ep(\om,\bk,\ep \xi) = \frac{i}{2 \sqrt{\beta(\bk)}
  \eta_s v_{oz}} \widecheck{S}\Big( -\frac{\xi}{ \eta_s v_{oz}},
\frac{\bk}{\gamma_s} \Big)e^{-i \frac{\om \xi}{ v_{oz}}- i \beta(\bk)
  \xi} + O(\sqrt{\ep}).
\end{align*}
Integrating in $\xi$ and using that $a^\ep(\om,\bk,\zeta)$ vanishes
for $\zeta \ll - O(\ep)$, we obtain that
\begin{align*}
a^\ep(\om,\bk,\ep \xi) &\approx \frac{i}{2  \sqrt{\beta(\bk)} \eta_s
  v_{oz}} \int_{\RR} {\rm d} \xi \, \widecheck{S}\Big( -\frac{\xi}{
  \eta_s v_{oz}}, \frac{\bk}{\gamma_s} \Big)e^{-i \frac{\om \xi}{
    v_{oz}}- i \beta(\bk) \xi} \nonumber \\ &= \frac{i}{2
  \sqrt{\beta(\bk)}} \hat S \Big(\eta_s(\om +
\beta(\bk)v_{oz}),\frac{\bk}{\gamma_s}\Big).
\end{align*}
We use this expression as the initial condition for the forward going 
amplitudes
\begin{equation}
a^\ep(\om,\bk,0+) \approx \frac{i}{2\sqrt{\beta(\bk)}} \hat S
\Big(\eta_s(\om + \beta(\bk)v_{oz}),\frac{\bk}{\gamma_s}\Big),
\label{eq:An35}
\end{equation}
and drop the source term and the backward going amplitudes $a_{-}^\ep$
in equation \eqref{eq:An30} for range $\zeta > 0$.

\subsection{The acoustic pressure field in the Markovian limit}
\label{sect:Anal6}
By definitions \eqref{eq:P1}, \eqref{eq:An1}, \eqref{eq:An20} and the
scaling relations \eqref{eq:An6}, the acoustic pressure is
\begin{equation*}
p(T_Lt,L\bx,Lz) \approx - \frac{ {\rho_o}}{T_L}\partial_t
u^\ep(t,\bx,z-\ep v_{oz} t).
\end{equation*}
Furthermore, equation \eqref{eq:An31f} and the Fourier decomposition
\eqref{eq:An23} give that
\begin{equation*}
\hspace{-0.15in}\frac{p(T_Lt,L \bx,Lz)}{2 \pi c_o  {\rho_o}/\la_o} \approx
\int_{\RR} \frac{\rm{d} \om}{2 \pi} \int_{\RR^d} \frac{{\rm d} \bk}{(2
  \pi \ep)^d} \frac{i a^\ep(\om,\bk,z)}{\sqrt{\beta(\bk)}} e^{ - i
  \big(\frac{\om_o}{\ep} + \om + \beta(\bk) v_{oz}\big)t+ i
  \frac{(\bk, \beta(\bk))}{\ep} \cdot (\bx,z)},
\end{equation*}
where we have used equation \eqref{eq:An30} to write
$
a^\ep(\om,\bk,z-\ep v_{oz}t) = a^\ep(\om,\bk,z) + O(\sqrt{\ep})$.

The shifted scaled frequency $\om + \beta(\bk) v_{oz}$ appears in the
initial condition \eqref{eq:An35}, and the random processes $\hat
\nu_c$ and $\hat Q^\ep$ in equation \eqref{eq:An30} depend on
$\bk/\gamma$. Thus, it is convenient to introduce the variables
\begin{equation}
\Omega = \eta \big[ \om + \beta(\bk) v_{oz}\big], \quad \bK = \bk/\gamma,
\label{eq:An38}
\end{equation}
and rewrite the expression of the pressure as
\begin{equation}
\hspace{-0.11in}\frac{p(T_Lt,L \bx,Lz)}{\om_o {\rho_o}} \approx
\int_{\RR} \frac{{\rm d} \Omega}{2 \pi \eta}
\int_{\RR^d} \hspace{-0.02in}\frac{{\rm d} \bK}{(2 \pi \ep/\gamma)^d}
\frac{A^\ep(\Omega,\bK,z)}{\sqrt{\beta(\gamma \bK)}} e^{ - i
  \big(\frac{\om_o}{\ep} + \frac{\Omega}{\eta} \big)t+ i \frac{(\gamma \bK,
    \beta(\gamma \bK))}{\ep} \cdot (\bx,z)} ,
\label{eq:An40}
\end{equation}
 with redefined amplitude
\begin{equation}
A^\ep(\Omega,\bK,z) = i  a^\ep  \Big(\frac{\Omega}{\eta} -
\beta(\gamma \bK)v_{oz},\gamma \bK,z\Big) + o(1).
\label{eq:An39}
\end{equation}
The $o(1)$ term, which tends to zero as $\ep \to 0$, is used in this
definition so that we have an equal sign in the evolution equation for
$A^\ep$, derived from \eqref{eq:An30}, after neglecting the backward
going amplitudes,
\begin{align}
& \partial_z A^\ep(\Omega,\bK,z) = \frac{i
    k_o}{\beta(\gamma \bK)} \Big[\frac{\Omega}{\eta
      c_o}-\frac{\vv_o}{c_o} \cdot \big(\gamma \bK, \beta(\gamma
    \bK)\big)\Big]A^\ep(\Omega,\bK,z)  \nonumber \\ 
    &
      +\frac{i}{2} \sqrt{\frac{\gamma}{\ep}} \int_{\RR} \frac{{\rm d}
    \Omega'}{2\pi} \int_{\RR^d} \frac{{\rm d} \bK'}{(2 \pi)^d}
  \frac{A^\ep(\Omega',\bK',z)}{\sqrt{\beta(\gamma \bK) \beta(\gamma
      \bK')}}e^{i \big[\beta(\gamma \bK')-\beta(\gamma \bK)\big]
    \frac{z}{\ep}}  \nonumber \\
    & 
    \quad \times \Big[k_o^2  \bar{\sigma}_c\hat
    \nu_c\Big(\Omega - \Omega' - \eta \big(\gamma \bK-\gamma
    \bK',\beta(\gamma \bK)-\beta(\gamma \bK'))\cdot \vv_o, \bK-\bK',
    \frac{\gamma z}{\ep} \Big) \nonumber
    \\ &\quad-\frac{\gamma^2\bar{\sigma}_\rho}{2} \hat
    Q^{\ep}\Big(\Omega - \Omega' - \eta \big(\gamma \bK-\gamma
    \bK',\beta(\gamma \bK)-\beta(\gamma \bK'))\cdot \vv_o, \bK-\bK',
    \frac{\gamma z}{\ep} \Big) \Big],
\label{eq:An41}
\end{align}
for $z > 0$.  The initial condition \eqref{eq:An35} becomes
\begin{equation}
A^\ep(\Omega,\bK,0+) = A_o(\Omega,\bK) = -\frac{1}{2
  \sqrt{\beta(\gamma \bK)}} \hat S \Big(\frac{\eta_s}{\eta} \Omega
,\frac{\gamma}{\gamma_s} \bK\Big).
\label{eq:An42}
\end{equation}

\subsection{The Markovian limit}
\label{sect:MarkLim}
Let $L^2(\mathscr{O},\mathbb{C})$ be the space of complex-valued,
square-integrable functions defined on the set 
\begin{equation}
\label{def:calO}
\mathscr{O} =\{\Omega \in
\RR\} \times \{\bK\in \RR^d, \gamma |\bK|< k_o\}
\end{equation} and denote ${\itbf
  A}^\ep (z) = (A^\eps(\Omega,\bK,z))_{(\Omega,\bK)\in
  \mathscr{O}}$ for $z \ge 0$. From equation \eqref{eq:An41} we obtain the
conservation of energy relation
\begin{equation}
\partial_{z} 
\int_{\mathscr{O}} {\rm d} \Omega  {\rm d} \bK \,
|\tae(\Omega,\bK,{z})|^2 = 0,
\end{equation}
so the Markov process ${\itbf A}^\ep (z) \in
L^2(\mathscr{O},\mathbb{C})$ lives on the surface of the ball with
center at the origin and $\ep$ independent radius $R_{{\itbf A}}$
defined by
\begin{equation}
R_{{\itbf A}}^2 =
\int_{\mathscr{O}} {\rm d} \Omega  {\rm d} \bK \,
 |\tae(\Omega,\bK,{z})|^2 =
 \int_{\mathscr{O}} {\rm d} \Omega  {\rm d} \bK \,
  |A_o(\Omega,\bK)|^2.
  \label{eq:radius}
\end{equation}

We describe the Markovian limit $\ep \to 0$ in Appendix
\ref{ap:DL}. The result is that the process of ${\itbf A}^\ep(z)$
converges weakly in ${\mathcal C}([0,\infty),L^2)$ to a Markov process
  whose infinitesimal generator can be identified.  The first and
  second moments of the limit process are described below.

\subsubsection{The mean amplitude}
\label{sect:mean}
The expectation of $\tae(\Omega,\bK,z)$ in the limit $\ep \to 0$ is
given by
\begin{equation}
\lim_{\ep \to 0} \EE \big[\tae(\Omega,\bK,z)\big] = A_o(\Omega,\bK)
e^{i \theta(\Omega,\bK) z + D(\bK) z },
\label{eq:mean}
\end{equation}
where
\begin{align}
\theta(\Omega,\bK) =\frac{k_o}{\beta(\gamma \bK)}
\Big[\frac{\Omega}{\eta c_o} - \frac{\vv_o}{c_o} \cdot (\gamma
  \bK,\beta(\gamma \bK)) \Big] + 
  \frac{\gamma^3 \bar{\sigma}_\rho^2}{8  \beta(\gamma \bK)}
\Delta_{\vr} \cR_{\rho \rho}(0,\vr)\big|_{\vr = {\bf 0}},
\label{eq:mean1}
\end{align}
is a real phase and
\begin{align}
D(\bK) =& - \int_{|\bK'| < k_o/\gamma}
  \frac{{\rm d}\bK'}{(2\pi)^d}\, \frac{1}{4
    \beta(\gamma\bK)\beta(\gamma\bK')} \int_{\RR^d} {\rm d} \br
  \int_0^\infty {\rm d}r_z \, e^{- i \big(\bK- \bK',\frac{\beta(\gamma
      \bK)-\beta(\gamma \bK')}{\gamma} \big)\cdot {\vr}}\nonumber
  \\ &\times \left\{k_o^4  \bar{\sigma}_c^2 {\cR}_{cc}(0,\vr) +
  \frac{\gamma^4 \bar{\sigma}_\rho^2}{4}
  \Delta_{\vr}^2\cR_{\rho\rho}(0,\vr)- {k_o^2 \gamma^2
    \bar{\sigma}_c \bar{\sigma}_\rho}\Delta_{\vr} \cR_{c\rho}(0,\vr)  
\right\},\label{eq:mean2}
\end{align}
with $\vr = (\br,r_z)$. 
 Moreover, $ \mbox{Re}
\big[ D(\bK)\big] < 0$, because
\[
\int_{\RR^{d+1}} {\rm d} \br   \, e^{- i  \vec{\bK} \cdot
      {\vr}}       \left\{k_o^4  \bar{\sigma}_c^2 {\cR}_{cc}(0,\vr) +
  \frac{\gamma^4 \bar{\sigma}_\rho^2}{4}
  \Delta_{\vr}^2\cR_{\rho\rho}(0,\vr)- {k_o^2 \gamma^2
     \bar{\sigma}_c\bar{\sigma}_\rho}\Delta_{\vr} \cR_{c\rho}(0,\vr)  
\right\} \ge 0,
\]
is the power spectral density of the process 
\ba\label{eq:P}
X(t,\vr)  =   k_0^2  \bar{\sigma}_c \nu_c(t,\vr)  -  
\frac{\bar{\sigma}_\rho  \gamma^2}{2} \Delta_{\vr}  \nu_\rho(t,\vr) , 
\ea
in the variable $\vr$ (for fixed $t$), which is non-negative by
Bochner's theorem. Thus, the mean amplitude decays on the range scale
\begin{equation}
\mathscr{S}(\bK) = - \frac{1}{\mbox{Re} \big[ D(\bK)\big]},
\label{eq:mean3}
\end{equation} 
called the scattering mean free path. 
In the relatively high frequency regime the damping is mainly
due to the fluctuations  of the wave speed, while in the relatively low frequency regime
the damping  is mainly
due to the fluctuations  of the density.
This damping is the mathematical
manifestation of the randomization of the wave due to cumulative
scattering. 

Recall that we have assumed   $\bar{\sigma}_\rho=O(1)$.  Therefore, 
  in the regime  $\gamma \ll 1$, 
the $ \cR_{cc}$ term dominates in \eqref{eq:mean2}.

\subsubsection{The mean intensity}
\label{sect:meanint}
The expectation of the intensity 
\begin{equation}
\label{eq:TR1}
I(\Omega,\bK,{z}) = \lim_{\ep \to 0} \EE \big[|\tae(\Omega,\bK,{z})|^2\big]
\end{equation}
satisfies
\begin{align}
\label{eq:TR2}
 \partial_{z} I(\Omega,\bK,{z}) = 
 \int_{\mathscr{O}} \frac{{\rm d} \Omega'}{2\pi} \frac{ {\rm d} \bK'}{(2\pi)^d} \,
  Q(\Omega,\Omega',\bK,\bK') \left[
   I(\Omega',\bK',{z})-I(\Omega,\bK,{z})\right],
\end{align}
for ${z} > 0$, with initial condition obtained from equation \eqref{eq:An42}:
\begin{equation}
I(\Omega,\bK,0) = |A_o(\Omega,\bK)|^2.
\label{eq:TR3}
\end{equation} 
Denoting  the power spectrum of $X(t, \vr)$  in (\ref{eq:P})  by
$P(\Omega, \vec{\bK})$, and letting 
$$
   \left( \tilde{\Omega}, \vec{\bK} \right) 
     =  \left(\Omega   - \eta
\big(\gamma \bK, \beta(\gamma \bK))\cdot
\vv_o
   ,\bK,\frac{\beta(\gamma \bK)}{\gamma}\right) ,
$$
the kernel in \eqref{eq:TR2} isgiven by 
$$
 Q(\Omega,\Omega',\bK,\bK')
 =  \frac{ P( \tilde{\Omega} -\tilde{\Omega}',  \vec{\bK} -  \vec{\bK}' ) }
 {4 \beta(\gamma\bK)\beta(\gamma \bK')} ,
$$
that is explicitly
  \begin{align}
Q&(\Omega,\Omega',\bK,\bK') = \left\{k_o^4  \bar{\sigma}_c^2 \widetilde
\cR_{cc} + \frac{\gamma^4 \bar\sigma_{\rho}^2}{4}\left[ |\bK-\bK'|^2 +
  \Big(\frac{\beta(\gamma \bK)-\beta(\gamma \bK')}{\gamma} \Big)^2
  \right]^2 \widetilde \cR_{\rho \rho}\right. \nonumber
\\&\left. +k_o^2 \gamma^2  \bar{\sigma}_c\bar{\sigma}_\rho \left[ |\bK-\bK'|^2 +
  \Big(\frac{\beta(\gamma \bK)-\beta(\gamma \bK')}{\gamma} \Big)^2
  \right] \widetilde \cR_{\rho c} \right\}\frac{1}
  {  4   \beta(\gamma\bK)\beta(\gamma \bK')} ,
\label{eq:IntKer}
\end{align}
where $\widetilde \cR_{cc}$ is the power spectral density
\eqref{eq:TR4}, evaluated as
\[
 \widetilde \cR_{cc} =  \widetilde \cR_{cc} \Big(\Omega - \Omega' - \eta
\big(\gamma \bK-\gamma \bK',\beta(\gamma \bK)-\beta(\gamma \bK'))\cdot
\vv_o, \bK-\bK', \frac{\beta(\gamma \bK)-\beta(\gamma \bK')}{\gamma}
\Big), 
\] 
and similar for $\widetilde \cR_{ \rho c}$ and $\widetilde \cR_{\rho
 \rho}$.   
 
Note  that  the kernel satisfies
\begin{equation}
\int_{\mathscr{O}} \frac{{\rm d} \Omega'}{2 \pi}  \frac{{\rm d}\bK'}{(2 \pi)^d}\,
  Q(\Omega,\Omega',\bK,\bK') = - 2
\mbox{Re}[\mathscr{Q}(\bK)] = \frac{2}{\mathscr{S}(\bK)},
\label{eq:TR5}
\end{equation}
where $\mathscr{S}(\bK)$ is the scattering mean free path defined in
\eqref{eq:mean3}. 

\subsubsection{The Wigner transform}
\label{sect:Wigner}
The wave amplitudes decorrelate at distinct frequencies $\Omega \ne
\Omega'$ and wave vectors $\bK \ne \bK'$, meaning that
\begin{equation}
\lim_{\ep \to 0} \EE \big[\tae(\Omega,\bK,{z})
  \overline{\tae}(\Omega',\bK',z)\big] = \lim_{\ep \to 0} \EE
\big[\tae(\Omega,\bK,{z})] \lim_{\ep \to 0}
\EE[\overline{\tae}(\Omega',\bK',z)\big].
\label{eq:Adecorrel}
\end{equation}
The right-hand side is the product of the means of the mode
amplitudes, which decay on the range scale defined by the scattering
mean free path \eqref{eq:mean3}.

However, the amplitudes are correlated for $|\Omega - \Omega'| =
O(\ep)$ and $|\bK-\bK'| = O(\ep)$.  We are interested in the second
moment
\[
\EE \Big[\tae\Big(\Omega,\bK + \frac{\ep \bq}{2},{z}\Big)
  \overline{\tae\Big(\Omega,\bK- \frac{\ep \bq}{2},z\Big)}\Big],
\]
whose Fourier transform in $\bq$ gives the energy density resolved
over frequencies and directions of propagation. This is  the Wigner
transform defined by
\begin{equation}
W^\ep(\Omega,\bK,\bx,{z}) = \int_{\mathbb{R}^d} \frac{{\rm d} \bq}{(2
  \pi)^d} e^{i \bq \cdot ( \nabla \beta(\gamma\bK) {z} + \bx)} \EE
\hspace{-0.03in}\left[ \tae\Big(\Omega,\bK + \frac{\ep
    \bq}{2},{z}\Big) \overline{\tae\Big(\Omega,\bK - \frac{\ep
      \bq}{2},{z}\Big)}\right].
\label{eq:defWigner}
\end{equation}
We show in Appendix  \ref{sec:appW} that the Wigner
transform converges in the limit $\ep \to 0$ to $W(\Omega,\bK,\bx,{z})$,
the solution of the transport equation
\begin{align}
\left[\partial_{z} - \nabla\beta(\gamma \bK) \cdot \nabla_{\bx}
  \right] W(\Omega,\bK,\bx,{z}) =&
\int_{\mathscr{O}} \frac{{\rm d} \Omega'}{2 \pi}  \frac{{\rm d}\bK'}{(2 \pi)^d}\,
Q(\Omega,\Omega',\bK,\bK') \nonumber \\ &\times \left[
  W(\Omega',\bK',\bx,{z})-W(\Omega,\bK,\bx,{z})\right],
\label{eq:TR6}
\end{align}
for ${z} > 0$, with initial condition
\begin{equation}
W(\Omega,\bK,\bx,0) = |A_o(\Omega,\bK)|^2 \delta(\bx)
\label{eq:TR7}
\end{equation}

The transport equation \eqref{eq:TRresult} in the physical scales is
obtained from \eqref{eq:TR7} as explained in Appendix
\ref{ap:unscaled}. 
In the next section we will show how this equation simplifies in 
the paraxial  regime, when $\gamma \ll 1$. This is the result used
for the imaging applications discussed in section \ref{sect:Inverse}.

\subsection{The paraxial limit}
\label{sect:Anal7}
Equation \eqref{eq:TR6} shows that the energy is transported on the
characteristic
\begin{equation}
\bx = - \gamma \frac{\bK}{\beta(\gamma \bK)} {z},
\label{eq:TR7p}
\end{equation}
parametrized by $z$, and depending on the wave-vector $\bK$. Here 
$|\bx|/{z} = O(\gamma)$ quantifies the opening angle of the cone
(beam) of propagation with axis $z$.  We write this explicitly as 
\begin{equation}
\bX = \bx/\gamma, \quad \mbox{where} ~ |\bX| = O(1).
\label{eq:ParX}
\end{equation}
The paraxial regime corresponds
to a narrow beam, modeled by $\gamma \to 0$ and
\begin{equation}
\Gamma= \gamma/\gamma_s = O(1).
\label{eq:TR7pp}
\end{equation}

At the range $z = 0$ of the source  we have from \eqref{eq:TR7} and
\eqref{eq:An42} that 
\begin{equation}
 W(\Omega,\bK,\bx,0) = \frac{\big|\hat S\big( \frac{\eta_s}{\eta}
   \Omega,\Gamma\bK\big)\big|^2}{4 \gamma^d k_o } \delta(\bX),
\label{eq:TR10pp}
\end{equation}
 and to obtain a finite limit as $\gamma \to 0$ we rescale the Wigner transform as 
\begin{equation}
{\cal W}(\Omega,\bK,\bX,{z}) = \gamma^d W(\Omega,\bK,\gamma \bX,{z}).
\label{eq:Wparax}
\end{equation}
We also change variables in \eqref{eq:TR6},
\begin{align*}
 \Omega - \Omega' - \eta \big(\gamma(\bK-\bK'),\beta(\gamma
 \bK)-\beta(\gamma\bK')\big)\cdot \vv_o \leadsto \Omega', \qquad
 \bK-\bK' \leadsto \bK', 
\end{align*}
 and obtain the transport equation 
 \begin{align}
&\hspace{-0.1in}\left[\partial_{z} +\frac{\bK}{\beta(\gamma \bK)}
    \cdot \nabla_{\bX} \right] {\cal W} (\Omega,\bK,\bX,{z}) =
  \frac{k_o^4}{4 }
  \int_{\mathscr{O}} \frac{{\rm d} \Omega'}{2 \pi}  \frac{{\rm d}\bK'}{(2 \pi)^d}\, 
  \frac{1}{
    \beta(\gamma \bK) \beta(\gamma \bK')} \nonumber\\ & \hspace{0.1in}
  \times \Big[ \bar{\sigma}_c^2 \widetilde \cR_{cc} \Big(\Omega',\bK',
    \frac{\beta(\gamma \bK) - \beta(\gamma \bK - \gamma \bK')}{\gamma}
    \Big) + O(\gamma^2)\Big]\nonumber \\ &\hspace{0.1in}\times \Big[
    {\cal W}\left( \Omega - \Omega' - \eta \gamma \bv_o \cdot \bK' -
    \eta v_{oz} (\beta(\gamma \bK) - \beta(\gamma \bK - \gamma \bK')),
    \bK-\bK',\bX,{z}\right)  \nonumber \\ &\hspace{3.4in}
    - {\cal W}(\Omega,\bK,\bX,{z})\Big],
\label{eq:TR11}
\end{align}
for ${z} > 0$ and a finite $\gamma \ll 1$, where $O(\gamma^2)$ denotes
the $\widetilde \cR_{\rho c}$ and $\widetilde \cR_{\rho \rho}$ terms
in the kernel \eqref{eq:IntKer}.

Recall that $\gamma \ll 1$ so,  in order to observe a significant 
 effect of the ambient motion, we  rescale the transversal speed as 
\begin{equation}
\bv_o = \frac{{\itbf V}_o}{\eta \gamma}, \quad \mbox{with} ~ |{\itbf
  V}_o| = O(1).
\label{eq:TR12}
\end{equation}
With a similar scaling of the range velocity
\begin{equation}
v_{oz} = \frac{V_{oz}}{\eta \gamma}, \quad \mbox{with} ~ |V_{oz}| =
O(1),
\label{eq:TR13}
\end{equation}
we obtain that the range motion plays no role in equation
\eqref{eq:TR11} as $\gamma \to 0$, because
\[
\beta(\gamma \bK) = k_o + O(\gamma^2), \qquad \beta(\gamma \bK) -
\beta(\gamma \bK - \gamma \bK') = O(\gamma^2).
\]

The transport equation satisfied by the Wigner transform ${\cal
  W}(\Omega,\bK,\bX,{z})$ in the paraxial limit  $ \gamma \to 0$ is
\begin{align}
\nonumber \Big[\partial_{z} + \frac{\bK}{k_o} \cdot \nabla_{\bX} \Big]
          {\cal W}(\Omega,\bK,\bX,{z}) = \frac{k_o^2}{4 }
          \int_{\mathbb{R}^d} \frac{{\rm d} \bK'}{(2 \pi)^d}
          \int_{\mathbb{R}}\frac{{\rm d} \Omega'}{2 \pi} \,  \bar{\sigma}_c^2\widetilde
          \cR_{cc} (\Omega',\bK', 0) \nonumber \\ \times \big[
            \mathcal{W}\left( \Omega - \Omega' - \bK' \cdot {\itbf
              V}_o, \bK-\bK',\bX,{z}\right)
            -\mathcal{W}(\Omega,\bK,\bX,{z})\big],
\label{eq:TR15}
\end{align}
for $z > 0$, with initial condition
\begin{equation}
\mathcal{W}(\Omega,\bK,\bX,0) = \frac{\big|\hat S\big(
  \frac{\eta_s}{\eta} \Omega,\Gamma\bK\big)\big|^2}{4k_o}
\delta(\bX).
\label{eq:TR16}
\end{equation}
The transport equation \eqref{eq:TRresultP} in the physical scales is
obtained from \eqref{eq:TR15} using the scaling relations explained in
Appendix \ref{ap:unscaled}.

\section{Summary}
\label{sect:sum}
We introduced an analysis of sound wave propagation in a time
dependent random medium which moves due to an ambient flow at speed
$\vv(t,\vx)$, and is modeled by the wave speed $c(t,\vx)$ and mass
density $\rho(t,\vx)$. The random fields $\vv(t,\vx)$, $c(t,\vx)$ and
$\rho(t,\vx)$ have small, statistically correlated fluctuations about
the constant values $\vv_o$, $c_o$ and $\rho_o$, on the length scale
$\ell$ and time scale $T$. The analysis starts from Pierce's equation,
which is obtained from the linearization of the fluid dynamics
equations about an ambient flow, and applies to waves with central
wavelength $\la_o \ll \ell$. The excitation is from a stationary
source with radius $\ell_s$, which emits a narrowband signal of
duration $T_s$.  

The analysis is in a forward wave propagation regime to a large
distance (range) $L \gg \ell$, within a cone with small opening
angle. Using the diffusion approximation theory, we showed that the
coherent part (the expectation) of the wave decays exponentially in
$L/{\mathscr S}$, and quantified the frequency- and wavevector-dependent 
scattering mean free path ${\mathscr S}$. We also derived transport equations for the energy
density (Wigner transform) of the wave, which show explicitly the
effect of the ambient flow and net scattering in the time dependent
random medium.

We used the wave propagation theory to study the inverse problem of
localizing (imaging) the source from measurements at a stationary
array of receivers located at range $L$.  This study is in the regime
of paraxial wave propagation, where the Wigner transform can be
computed explicitly, and assumes a large range $L \gg {\mathscr S}$,
so that the wave is incoherent due to strong scattering in the random
medium.  The temporal variation of the medium is at time scale $T \ll
T_s$, and it has two beneficial effects for imaging: First, it causes
broadening of the bandwidth of the recorded waves, which leads to
improved travel time estimation and consequently, better range
resolution. Second, it allows a robust (statistically stable)
estimation of the Wigner transform from the array measurements.  We
presented an explicit analysis of imaging based on this Wigner
transform and showed how one can estimate the source location, the
mean velocity $\vv_o$ and the statistics of the random medium.

\section*{Acknowledgments}
Liliana Borcea's research  is supported in part by the Air Force Office of
Scientific Research under award number FA9550-18-1-0131 and in part by
the U.S. Office of Naval Research under award number N00014-17-1-2057.
Knut S\o lna's  research  is supported in part by the Air Force Office of
Scientific Research under award number FA9550-18-1-0217 and 
NSF grant 1616954. 
\appendix
\section{The Markovian limit theorem}
\label{ap:DL}
In this appendix we obtain the $\ep \to 0$ limit of the Markov process
${\itbf A}^\ep (z) = (A^\eps(\Omega,\bK,z))_{(\Omega,\bK)\in
  \mathscr{O}}$, which lies on the surface of the sphere with radius
$R_{\itbf A}$ given in equation \eqref{eq:radius}. The set
$\mathscr{O}$ is defined by (\ref{def:calO}). 
The process ${\itbf A}^\ep (z)$ starts from
\begin{equation}
{\itbf A}^\ep (0) = \left(-\frac{1}{2 \sqrt{\beta(\gamma\bK})} \hat S
\Big(\frac{\eta_s}{\eta} \Omega ,\frac{\gamma}{\gamma_s} \bK\Big)
\right)_{(\Omega,\bK)\in \mathscr{O}},
\label{eq:startA}
\end{equation}
which is independent of $\ep$, and evolves at $z >0$ according to the
stochastic equation
\begin{equation}
\frac{{\rm d} {\itbf A}^\ep}{{\rm d}z} = {\mathcal
  G}\Big(\frac{z}{\ep},\frac{z}{\ep}\Big) {\itbf A}^\ep +
\frac{1}{\sqrt{\ep}} {\mathcal F}\Big(\frac{z}{\ep},\frac{z}{\ep}\Big)
     {\itbf A}^\ep.
\label{eq:StochA}
\end{equation}
Here ${\mathcal G}$ and ${\mathcal F}$ are integral operators
\begin{align*}
[{\mathcal G}(z,\zeta) {\itbf A}](\Omega,\bK)&= \int_{\mathscr{O}}
{\rm d}\Omega' {\rm d}\bK'\, G(z,\zeta,\Omega,\Omega',\bK,\bK')
A(\Omega',\bK') , \\ [{ \mathcal F}(z,\zeta) {\itbf A}](\Omega,\bK)&=
\int_{\mathscr{O}}{\rm d}\Omega' {\rm d}\bK'\,
F(z,\zeta,\Omega,\Omega',\bK,\bK') A(\Omega',\bK'),
\end{align*}
with kernels depending on the random processes $\nu_c(\tau,\vr)$ and
$Q^\ep(\tau,\vr)$. Recall the definition \eqref{eq:An16} of $Q^\ep(\tau,\vr)$. We rewrite it
here as
\begin{equation}
Q^\ep(\tau,\vr) = Q^{(0)}(\tau,\vr) + \frac{\sqrt{\ep \gamma}
  \bar{\sigma}_\rho}{2}Q^{(1)}(\tau,\vr),  
\label{eq:QepDef}
\end{equation}
with 
\begin{equation}
Q^{(0)}(\tau,\vr) = \Delta_{\vr} \nu_\rho(\tau,\vr) ~~ \mbox{and} ~~
Q^{(1)}(\tau,\vr) = \left|\nabla_{\vr} \nu_\rho(\tau,\vr) \right|^2.
\label{eq:Q01Def}
\end{equation}
The kernel $G$ has a deterministic part supported at $\Omega' =
\Omega$ and $\bK' = \bK$, and a random part determined by $Q^{(1)}$,
\begin{align}
G(z,\zeta,&\Omega,\Omega',\bK,\bK') = \frac{i k_o}{\beta(\gamma \bK)}
\Big[\frac{\Omega}{\eta c_o}-\frac{\vv_o}{c_o} \cdot \big(\gamma \bK,
  \beta(\gamma \bK)\big)\Big] \delta(\Omega'-\Omega) \delta(\bK'-\bK)
\nonumber \\ &-\frac{i \gamma^3 \bar{\sigma}_\rho^2}{8(2
  \pi)^{d+1}\sqrt{\beta(\gamma \bK) \beta(\gamma \bK')}}
\exp\big\{ i \big[ \beta(\gamma \bK') - \beta(\gamma
    \bK)\big]\zeta\big\} \nonumber \\ &\times\hat
Q^{(1)}\Big(\Omega - \Omega' - \eta \big(\gamma \bK-\gamma
\bK',\beta(\gamma \bK)-\beta(\gamma \bK'))\cdot \vv_o, \bK-\bK',
\gamma z \Big). \label{eq:DefkerG}
 \end{align}
The kernel $F$ is determined by $\nu_c$ and $Q^{(0)}$, 
\begin{align}
\hspace{-0.1in}F(&z,\zeta,\Omega,\Omega',\bK,\bK') = \frac{i
  \sqrt{\gamma}}{2(2\pi)^{d+1}\sqrt{\beta(\gamma \bK)\beta( \gamma
    \bK')}} \exp\big\{ i [\beta(\gamma \bK')-\beta(\gamma\bK)] \zeta
\big\} \nonumber \\ &\hspace{-0.07in}\times\left[ k_o^2
  \bar{\sigma}_c  \hat{\nu}_c\Big(\Omega - \Omega' - \eta \big(\gamma \bK-\gamma
  \bK',\beta(\gamma \bK)-\beta(\gamma \bK'))\cdot \vv_o, \bK-\bK',
  \gamma z \Big) \right. \nonumber
  \\ &\hspace{-0.07in}\left. -\frac{\gamma^2 \bar{\sigma}_\rho}{2} \hat
  Q^{(0)}\Big(\Omega - \Omega' - \eta \big(\gamma \bK-\gamma
  \bK',\beta(\gamma \bK)-\beta(\gamma \bK'))\cdot \vv_o, \bK-\bK',
  \gamma z \Big)\right].\label{eq:DefkerF}
\end{align}

The random process ${\itbf A}^\eps(z)$ is Markov with generator
\begin{align*}
{\cal L}^\eps f({\itbf A},\overline{{\itbf A}}) =&
\int_{\mathscr{O}^2} \frac{1}{\sqrt{\ep}}
F\Big(\frac{z}{\ep},\frac{z}{\ep},\Omega,\Omega',\bK,\bK' \Big)
\frac{\delta f}{\delta {A}(\Omega,\bK)} {A} (\Omega',\bK') {\rm
  d}\Omega' {\rm d}\bK' {\rm d}\Omega {\rm d}\bK \\ & +
\int_{\mathscr{O}^2} \frac{1}{\sqrt{\ep}}
\overline{F}\Big(\frac{z}{\ep},\frac{z}{\ep},\Omega,\Omega',\bK,\bK'
\Big) \frac{\delta f}{\delta \overline{{A}}(\Omega,\bK)}
\overline{{A}} (\Omega',\bK') {\rm d}\Omega' {\rm d}\bK'{\rm d}\Omega
         {\rm d}\bK\\ &+ \int_{\mathscr{O}^2}
         G\Big(\frac{z}{\ep},\frac{z}{\ep},\Omega,\Omega',\bK,\bK'
         \Big) \frac{\delta f}{\delta {A}(\Omega,\bK)} {A}
         (\Omega',\bK') {\rm d}\Omega' {\rm d}\bK' {\rm d}\Omega {\rm
           d}\bK \\ & + \int_{\mathscr{O}^2}
         \overline{G}\Big(\frac{z}{\ep},\frac{z}{\ep},\Omega,\Omega',\bK,\bK'
         \Big) \frac{\delta f}{\delta \overline{{A}}(\Omega,\bK)}
         \overline{{A}} (\Omega',\bK') {\rm d}\Omega' {\rm d}\bK'{\rm
           d}\Omega {\rm d}\bK ,
\end{align*}
where $\delta f/\delta {A}(\Omega,\bK)$ denotes the variational
derivative, defined as follows.  If $\varphi$ is a smooth function and
\begin{align*}
f({\itbf A},\overline{{\itbf A}}) = 
\int \cdots \int \varphi(\Omega_1,\ldots,\Omega_{n+m},
\bK_1,\ldots, \bK_{n+m}) 
\prod_{j=1}^n {A}(\Omega_j,\bK_j) \\
\times \prod_{j=n+1}^{n+m} \overline{{A}}(\Omega_{j},\bK_{j})
\prod_{j=1}^{n+m} {\rm d}\Omega_j {\rm d}\bK_j  ,
\end{align*}
then  we have
\begin{align*}
\frac{\delta f}{\delta {A}(\Omega,\bK)} = 
\sum_{l=1}^n
\int \cdots \int \varphi(\Omega_1,\ldots,\Omega_{n+m},
\bK_1,\ldots, \bK_{n+m}) \mid_{_{\Omega_l=\Omega,\bK_l=\bK}}
\hspace{-0.05in} \prod_{j=1,j\neq l}^n \hspace{-0.05in}
       {A}(\Omega_j,\bK_j) \\ \times
       \prod_{j=n+1}^{n+m} \hspace{-0.05in}
       \overline{{A}}(\Omega_{j},\bK_{j}) \prod_{j=1, j\neq
         l}^{n+m} \hspace{-0.05in} {\rm d}\Omega_j {\rm d}\bK_j
\end{align*}
and
\begin{align*}
\frac{\delta f}{\delta \overline{{A}}(\Omega,\bK)} = 
\sum_{l=n+1}^{n+m}
\int \cdots \int \varphi(\Omega_1,\ldots,\Omega_{n+m},
\bK_1,\ldots, \bK_{n+m}) \mid_{_{\Omega_l=\Omega,\bK_l=\bK}}
\prod_{j=1}^n {A}(\Omega_j,\bK_j) \\ \times
\prod_{j=n+1,j\neq l}^{n+m} \overline{{A}}(\Omega_{j},\bK_{j})
\prod_{j=1, j\neq l}^{n+m} {\rm d}\Omega_j {\rm d}\bK_j .
\end{align*}
The linear combinations of such functions $f$ form an algebra that is
dense in
${\cal C}(L^2)$ 
 and is convergence determining.  We can also
extend the class of functions to include generalized functions
$\varphi$ of the form
\begin{align*} 
\varphi(\Omega_1,\ldots,\Omega_{2n},\bK_1,\ldots,\bK_{2n}) = &
\Phi(\Omega_1,\ldots,\Omega_n,\bK_1,\ldots,\bK_n) \\
& \times \prod_{j=1}^n
\delta(\Omega_{n+j}-\Omega_j) \delta(\bK_{n+j}-\bK_j) ,
\end{align*}
where $\Phi$ is a smooth function.

Applying the diffusion-approximation theory described in \cite[Chapter
  6]{book} and \cite{PW1994,PK1974}, we obtain the limit generator
\begin{align}
{\cal L} &f({\itbf A},\overline{{\itbf A}}) = \int_0^\infty {\rm
  d}\zeta \lim_{Z\to \infty}\frac{1}{Z} \int_0^Z {\rm d}h
\int_{\mathscr{O}^4}{\rm d} \Omega_1'{\rm d}\bK_1' {\rm
  d}\Omega_2'{\rm d}\bK_2' {\rm d} \Omega_1{\rm d}\bK_1 {\rm
  d}\Omega_2{\rm d}\bK_2\nonumber\\ &\times \Big\{ \EE\big[
  F(0,h,\Omega_1,\Omega_1',\bK_1,\bK_1')
  F(\zeta,\zeta+h,\Omega_2,\Omega_2',\bK_2,\bK_2') \big]
\nonumber\\ &\quad \quad \times \frac{\delta^2 f}{\delta
  {A}(\Omega_1,\bK_1)\delta {A}(\Omega_2,\bK_2)}
 {A}(\Omega_1',\bK_1'){A}(\Omega_2',\bK_2')
\nonumber\\ &\quad+ \EE\big[
F(0,h,\Omega_1,\Omega_1',\bK_1,\bK_1')
\overline{F}(\zeta,\zeta+h,\Omega_2,\Omega_2',\bK_2,\bK_2')
\big] \nonumber\\&\quad \quad \times\frac{\delta^2
f}{\delta {A}(\Omega_1,\bK_1)\delta
\overline{{A}}(\Omega_2,\bK_2)}
{A}(\Omega_1',\bK_1')\overline{{A}}(\Omega_2',\bK_2')\nonumber\\ &
\quad+ \EE\big[\overline{F}(0,h,\Omega_1,\Omega_1',\bK_1,\bK_1')
 {F}(\zeta,\zeta+h,\Omega_2,\Omega_2',\bK_2,\bK_2')\big] \nonumber\\&\quad
\quad\times\frac{\delta^2 f}{\delta \overline{{A}}(\Omega_1,\bK_1)\delta
{{A}}(\Omega_2,\bK_2)} \overline{{A}}(\Omega_1',\bK_1')
{{A}}(\Omega_2',\bK_2')\nonumber\\ &\quad
 + \EE\big[\overline{F}(0,h,\Omega_1,\Omega_1',\bK_1,\bK_1') 
\overline{F}(\zeta,\zeta+h,\Omega_2,\Omega_2',\bK_2,\bK_2')
 \big] \nonumber \\&\quad \quad\times\frac{\delta^2
f}{\delta\overline{{A}}(\Omega_1,\bK_1)\delta
 \overline{{A}}(\Omega_2,\bK_2)} \overline{{A}}(\Omega_1',\bK_1')
\overline{{A}}(\Omega_2',\bK_2') \Big\} \nonumber\\ & + \int_0^\infty
{\rm d}\zeta \lim_{Z\to \infty}\frac{1}{Z} \int_0^Z {\rm d}h \int_{\mathscr{O}^3}{\rm d}
 \Omega_1'{\rm d}\bK_1'{\rm d} \Omega_1{\rm d}\bK_1 {\rm d}\Omega'{\rm d}\bK'
\nonumber\\ &\times \Big\{ \EE\big[ {F}(0,h,\Omega',\Omega_1',\bK',\bK_1')
{F}(\zeta,\zeta+h,\Omega_1,\Omega',\bK_1,\bK')
\big]
\frac{  \delta  f}{\delta
{{A}}(\Omega_1,\bK_1)}{{A}}(\Omega_1',\bK_1') \nonumber\\ & \quad 
+ \EE\big[
\overline{F}(0,h,\Omega',\Omega_1',\bK',\bK_1')
 \overline{F}(\zeta,\zeta+h,\Omega_1,\Omega',\bK_1,\bK') \big]
 \frac{\delta f}{\delta \overline{{A}}(\Omega_1,\bK_1)}
\overline{{A}}(\Omega_1',\bK_1') \Big\}\nonumber \\ & + \lim_{Z\to
  \infty}\frac{1}{Z} \int_0^Z {\rm d}h \int_{\mathscr{O}^2}{\rm d} \Omega{\rm d}\bK{\rm d}
\Omega'{\rm d}\bK' \Big\{ \EE\big[ {G}(0,h,\Omega,\Omega',\bK,\bK')\big]
 \nonumber \\ & \quad \times \frac{\delta f}{\delta {{A}}(\Omega,\bK)}
 {{A}}(\Omega',\bK')+ \EE\big[ \overline{G}(0,h,\Omega,\Omega',\bK,\bK') \big]
  \frac{\delta  f}{\delta \overline{{A}}(\Omega,\bK)} \overline{{A}}(\Omega',\bK')
 \Big\}. \label{eq:Generator}
 \end{align}
The expectations in the expression of the generator can be computed
with
\begin{align}
\EE\big[ \hat{\nu}_c(\Omega,\bK,\zeta)\hat{\nu}_c(\Omega',\bK',0)\big]
=& (2\pi)^{d+1} \delta(\bK+\bK') \delta(\Omega+\Omega') \hat{\cR}_{cc}
(\Omega,\bK,\zeta) , \label{eq:nunu}\\ \EE\big[
  \overline{\hat{\nu}_c}(\Omega,\bK,\zeta)\hat{\nu}_c(\Omega',\bK',0)\big]
=& (2\pi)^{d+1} \delta(\bK-\bK') \delta(\Omega-\Omega')
\overline{\hat{\cR}_{cc} (\Omega,\bK,\zeta)} , \\ \EE\big[
  \hat{\nu}_c(\Omega,\bK,\zeta)\overline{\hat{\nu}_c}(\Omega',\bK',0)\big]
=& (2\pi)^{d+1} \delta(\bK-\bK') \delta(\Omega-\Omega') \hat{\cR}_{cc}
(\Omega,\bK,\zeta) , \\ \EE\big[
  \overline{\hat{\nu}_c}(\Omega,\bK,\zeta)
  \overline{\hat{\nu}_c}(\Omega',\bK',0)\big] =& (2\pi)^{d+1}
\delta(\bK+\bK') \delta(\Omega+\Omega') \overline{\hat{\cR}_{cc}
  (\Omega,\bK,\zeta)}, \label{eq:nunu1}
\end{align}
and similar for $\hat \nu_\rho$. 
Note here that both  $\nu_c$ and $\hat{\cR}_{cc}$ are real.
We also have
\begin{align}
\hat Q^{(0)}(\Omega,\bK,z) = \big(-|\bK|^2 +\partial_z^2\big) \hat
\nu_\rho(\Omega,\bK,z),\label{eq:defHatQ0}
\end{align}
and 
\begin{align}
\EE\big[ \hat Q^{(1)}(\Omega,\bK,z)\big] = - (2 \pi)^{d+1}
\delta(\Omega) \delta(\bK) \Delta_{\vr} \cR_{\rho
  \rho}(0,\vr)\big|_{\vr = {\bf 0}}.
\label{eq:defMeanQ1}
\end{align}

\subsection{The mean amplitude}
To calculate the mean of the limit process, we let
\[f({\itbf A},\overline{\itbf A})=
\int_{\mathscr{O}} {\rm d}\Omega {\rm d}\bK\,
\varphi(\Omega,\bK){A}(\Omega,\bK),\] so that
\[\frac{\delta f}{\delta {{A}}(\Omega_1,\bK_1)} = \varphi(\Omega_1,\bK_1), \quad  \quad 
\frac{\delta f}{\delta \overline{{A}}(\Omega_1,\bK_1)} =0,\] and all
second variational derivatives are zero.  

From the expression \eqref{eq:Generator}, definitions
\eqref{eq:DefkerG}--\eqref{eq:DefkerF} and the expectations
\eqref{eq:nunu}--\eqref{eq:defMeanQ1} we obtain
\begin{align}
{\cal L} f({\itbf A}, \overline{\itbf A})&= \int_{\mathscr{O}} {\rm
  d}\Omega {\rm d}\bK \, \big[i \theta(\Omega,\bK)+
  D(\bK)\big]\varphi(\Omega,\bK) A(\Omega,\bK) ,
\end{align}
with $\theta$ and $D$ given in \eqref{eq:mean1} and \eqref{eq:mean2}. 
This gives the result \eqref{eq:mean}.

\subsection{The mean intensity}
To characterize the mean intensity of the limit process, we let
\begin{align*}
 f({\itbf A},\overline{\itbf A})& = \int_{\mathscr{O}} {\rm d}\Omega
 {\rm d}\bK \, \varphi(\Omega,\bK)|{A}(\Omega,\bK)|^2   
  \\&= \int_{\mathscr{O}^2}{\rm d}\Omega {\rm d}\bK{\rm
   d}\Omega' {\rm d}\bK'\, \varphi(\Omega,\bK) \delta(\Omega-\Omega')
 \delta(\bK-\bK'){A}(\Omega,\bK) \overline{A}(\Omega',\bK') ,
\end{align*}
so that 
\begin{align*}
& \frac{\delta f}{\delta
    {{A}}(\Omega_1,\bK_1)} = \overline{{A}}(\Omega_1,\bK_1)
  \varphi(\Omega_1,\bK_1),\quad \quad \frac{\delta f}{\delta
    \overline{{A}}(\Omega_1,\bK_1)} = {{A}}(\Omega_1,\bK_1)
  \varphi(\Omega_1,\bK_1),\\ 
  &\frac{\delta^2 f}{\delta
    \overline{{A}}(\Omega_1,\bK_1)\delta {{A}}(\Omega_2,\bK_2)} =
  \varphi(\Omega_2,\bK_2) \delta(\Omega_2-\Omega_1)
  \delta(\bK_2-\bK_1),
\end{align*} 
and all other second variational derivatives are zero.

Using the expectation \eqref{eq:defMeanQ1} and definition
\eqref{eq:DefkerG} in \eqref{eq:Generator}, we obtain that the $G$
dependent terms make no contribution. Furthermore, using the
expectations \eqref{eq:nunu}--\eqref{eq:nunu1} and \eqref{eq:defHatQ0}
we get
\begin{align*}
{\cal L}f({\itbf A},\overline{\itbf A})=& -\int_{\mathscr{O}}
\frac{{\rm d} \Omega_1{\rm d}\bK_1}{(2 \pi)^{d+1}} \,
\varphi(\Omega_1,\bK_1) |{{A}}(\Omega_1,\bK_1)|^2 \int_{\mathscr{O}} {\rm
  d}\Omega_1'{\rm d}\bK_1'
Q(\Omega_1,\Omega_1',\bK_1,\bK_1') \\ &+ \int_{\mathscr{O}}
\frac{{\rm d} \Omega_1{\rm d}\bK_1}{(2 \pi)^{d+1}} \,
\varphi(\Omega_1,\bK_1) \int_{\mathscr{O}} {\rm d}\Omega_1'{\rm d}\bK_1'
\, |{{A}}(\Omega_1',\bK_1')|^2
Q(\Omega_1,\Omega_1',\bK_1,\bK_1'),
\end{align*}
with kernel defined in \eqref{eq:IntKer}.
This gives the equation satisfied by the mean intensity.

\subsection{Wave decorrelation and the Wigner transform}\label{sec:appW}
To study the second moments at distinct frequencies $\Omega$,
$\Omega'$ and wave vectors $\bK$ and $\bK'$, we let
\begin{align*}
f({\itbf A},\overline{\itbf A})& = \int_{\mathscr{O}^2}{\rm d}\Omega {\rm d}\bK{\rm d}\Omega' {\rm
  d}\bK'\, 
\varphi(\Omega,\Omega',\bK,\bK'){A}(\Omega,\bK)
\overline{A}(\Omega',\bK') .
\end{align*}
Then, we have 
\begin{align*} \frac{\delta f}{\delta {{A}}(\Omega_1,\bK_1)} &=
\int_{\mathscr{O}} {\rm d}\Omega' {\rm d} \bK'\,
\overline{{A}}(\Omega',\bK') \varphi(\Omega_1,\Omega',\bK_1,\bK')
,\\ \frac{\delta f}{\delta \overline{{A}}(\Omega_1,\bK_1)} &=
\int_{\mathscr{O}} {\rm d} \Omega'{\rm d} \bK' \, {A}(\Omega',\bK')
\varphi(\Omega',\Omega_1,\bK',\bK_1), \\
\frac{\delta^2 f}{\delta{{A}}(\Omega_1,\bK_1)\delta  
\overline{{A}}(\Omega_2,\bK_2)} &=
\varphi(\Omega_1,\Omega_2,\bK_1,\bK_2) ,
\end{align*} 
and all other second variational derivatives are zero.

Substituting in \eqref{eq:Generator} and using the expectations
\eqref{eq:nunu}--\eqref{eq:defMeanQ1} we obtain that 
\begin{align}
{\cal L} f({\itbf A}, \overline{\itbf A})=& \int_{\mathscr{O}} {\rm
  d}\Omega {\rm d}\bK \int_{|\bK'| < k_o} {\rm d} \bK' \,
\varphi(\Omega,\Omega',\bK,\bK') \nonumber \\ &\times \big[i
  \theta(\Omega,\bK)- i \theta(\Omega',\bK') + D(\bK) +
  \overline{D}(\bK')\big] A(\Omega,\bK)\overline{A}(\Omega',\bK'),
\label{eq:decorrel}
\end{align}
with $\theta(\Omega,\bK)$ and $D(\bK)$ defined in 
\eqref{eq:mean1} and \eqref{eq:mean2}. 
This gives the decorrelation result \eqref{eq:Adecorrel}. 

Finally, to study the Wigner transform, we use \eqref{eq:An41} to
obtain an evolution equation for
\begin{align*}
\mathscr{W}^\ep(\Omega,\bK,t,\bx,z) =& \int_{\RR} \frac{{\rm d} w}{2
  \pi} \int_{\RR} \frac{{\rm d} \bq}{(2 \pi)^d} e^{-i w t + i \bq \cdot[
    \bx + \nabla \beta(\bK) z]} \\
 &\times    A^\ep\Big(\Omega + \frac{\ep w}{2},\bK
+ \frac{\ep \bq}{2}, z \Big) \overline{A}^\ep\Big(\Omega - \frac{\ep
  w}{2},\bK - \frac{\ep \bq}{2}, z \Big),  
 \end{align*}
and then analyze the limit $\ep \to 0$ of $\mathscr{W}^\ep$ with the
same approach as described in this appendix. The Wigner transform
\eqref{eq:defWigner} is
\[
W^\ep(\Omega,\bK,\bx,z) = \int_{\RR} {\rm d}t \, \EE
\big[\mathscr{W}^\ep(\Omega,\bK,t,\bx,z)\big],
\]
and this converges in the limit to the solution $W(\Omega,\bK,\bx,z)$
of \eqref{eq:TR6}--\eqref{eq:TR7}.

\section{The transport equation in the physical scales}
\label{ap:unscaled}

To distinguish between the scaled and unscaled variables, we resurrect
the notation of section \ref{sect:Anal2} with the  unscaled variables
denoted by primes.

We begin with the pressure field \eqref{eq:An40},
\begin{align}
p(T_Lt',L \bx',Lz') \approx \frac{i \om_o {\rho_o} }{\eps^d}
\int\frac{{\rm d} \omega'}{2\pi}\int \frac{{\rm d} \bk'}{(2\pi)^d}
\frac{ a^{\ep \,'}(\om' - \beta'(\bk') v'_z,\bk',z')
}{\sqrt{\beta'(\bk')}} \nonumber \\ \times e^{-i
  \big(\frac{\om'_o}{\ep} + \om' \big)t' +i \frac{\bk'}{\ep} \cdot
  \bx' + i\frac{\beta'(\bk')}{\eps} z' }.
    \label{eq:B1}
\end{align}
The scaling relations \eqref{eq:F8}--\eqref{eq:F13}
and \eqref{eq:An6}--\eqref{eq:An7} give
\begin{align*}
\frac{\om_o' t'}{\ep} &= \frac{2 \pi}{\la_o/L} \frac{t}{T_L} = \om_o
t,\\ \om' t' & = \om' \frac{t}{T_L} = \om t, \quad \mbox{i.e.,} ~ \om' =
\om T_L,\\ \frac{\bk'}{\ep} \cdot \bx' &= \frac{\bk'}{\la_o/L} \cdot
\frac{\bx}{L} = \frac{\bk'}{\la_o} \cdot \bx = \bk \cdot \bx, \quad
\mbox{i.e.,} ~ \bk' = \la_o \bk,\\ \beta'(\bk') &= \sqrt{(k_o')^2 -
  |\bk'|^2} = \la_o \sqrt{k_o^2 - |\bk|^2} = \la_o \beta(\bk),
\\ \om'-\beta'(\bk') v'_{oz} &= T_L \om - \la_o \beta(\bk)
  \frac{v_{oz}}{(\la_o/L) c_o}  = T_L [\om - \beta(\bk) v_{oz}],
\\ \frac{\beta'(\bk')}{\ep} z' &= \frac{\la_o \beta(\bk)}{\la_o/L}
\frac{z}{L} = \beta(\bk) z.
\end{align*}
Equation \eqref{eq:B1} becomes \eqref{eq:B3}, with amplitudes
\begin{equation}
a(\om,\bk,z) = \frac{T_L L^d}{\sqrt{\la_o}} a^{\ep \,'}
\left( {T_L}(\om-\beta(\bk)v_{oz}),\la_o \bk,\frac{z}{L} \right),
\label{eq:B4}
\end{equation}
satisfying the initial conditions
\begin{equation}
a(\om,\bk,0) = \frac{T_L L^d}{\la_o} \frac{i}{2 \sqrt{\beta(\bk)}}
\hat S(T_s \om, \ell_s \bk) = \frac{i\sigma_s T_s \ell_s^d}{2
  \sqrt{\beta(\bk)}}\hat S(T_s \om, \ell_s \bk),
\label{eq:B5}
\end{equation}
derived from \eqref{eq:F14} and  \eqref{eq:An35}.  

It remains to write the transport equation \eqref{eq:TRresult} for the
Wigner transform. To do so, we obtain from definitions \eqref{eq:An38}
and the scaling relations above that
\begin{align*}
\Omega' & = \eta [ \om' + \beta'(\bk') v'_{oz}] = {T} [ \om +
  \beta(\bk) v_{oz}] , \\ \bK'&= \frac{\bk'}{\gamma} = \frac{\la_o
  \bk}{\la_o/\ell} = \ell \bk.
\end{align*}
We also recall the definition \eqref{eq:An39} of $A^\ep$ in terms of
$a^\ep$, and obtain that
\begin{align*}
W(\om,&\bk,\bx,z) = \int \frac{{\rm d} \bq}{(2\pi)^d} \, \exp \Big[ i\bq
  \cdot \big( \nabla \beta(\bk) z + \bx \big)\Big] \EE
\left[a\Big(\om,\bk + \frac{\bq}{2}\Big) \overline{a\Big(\om,\bk -
    \frac{\bq}{2}\Big)} \right] \nonumber \\ 
    &=\left(\frac{T_L
  L^d}{\sqrt{\la_o}} \right)^2 \frac{1}{\la_o^d} \int \frac{{\rm d} \bq'}{(2
  \pi)^d} \exp \Big[ i\frac{\bq'}{\la_o/L} \cdot \big( \nabla
  \beta'(\bk') z' + \bx' \big)\Big] \\ &\times \EE
\left[A^\ep\Big(\Omega' -  \eta \beta'(\gamma\bK')v_{oz}',\bK' +
  \frac{\bq'}{2},z'\Big) \overline{A^\ep\Big(\Omega' - \eta
    \beta(\gamma\bK')v_{oz}',\bK' - \frac{\bq'}{2},z'\Big)} \right],
\end{align*}
with $\Omega'$ and $\bK'$ defined as above in terms of $\om$ and
$\bk$, and $z' = z/L$.  Since $\ep = \la_o/L$, we can change the
variable of integration as $\bq' \leadsto \bq' \ep$, and obtain
\begin{align}
W(\om,\bk,\bx,z) &= \frac{T_L^2 L^d}{\la_o} W^{\ep \, '}\big( \Omega'-
\eta \beta'(\gamma \bK') v_{oz}',\bK',\bx',z'\big) \nonumber
\\ &\approx  \sigma_s^2 T_s^2 \ell_s^{2d} \frac{\la_o}{L^d} W'\big(
\Omega'- \eta \beta'(\gamma \bK') v_{oz}',\bK',\bx',z'\big).
\end{align}
Here the approximation is for $\ep \ll 1$, where we have replaced $W^{\ep
  \, '}$ by its $\ep \to 0$ limit~$W^{'}$.

Using the initial conditions \eqref{eq:An42} and \eqref{eq:TR7} and the
scaling relations between $\Omega'$, $\bK'$ and $\om$ and $\bk$, we
have
\begin{align*}
W(\om,\bk,\bx,0) &= \sigma_s^2 T_s^2 \ell_s^{2d} \frac{\la_o}{L^d}
\frac{\delta(\bx/L)}{4 \beta'(\bk')} \Big|\hat
S\Big(\eta_s \om',\frac{\ell_s}{\ell}\bK'\Big)
\Big|^2 \\&= \sigma_s^2 T_s^2 \ell_s^{2d} \frac{\delta(\bx)}{4
  \beta(\bk)} \big|\hat S\big( T_s \om,\ell_s \bk\big) \big|^2,
\end{align*}
as stated in \eqref{eq:TRinitial} and (\ref{eq:R2}). The transport equation
\eqref{eq:TRresult} follows from \eqref{eq:TR6}, using
\[
\partial_{z'} - \nabla \beta'(\gamma \bK') \cdot \nabla_{\bx'} = L
\left[ \partial_z - \nabla \beta(\bk) \cdot \nabla _\bx\right].
\]

\section{Solution of the transport equation in the paraxial regime}
\label{ap:SolPar}
To deal with the convolution in \eqref{eq:TRresultP}, we Fourier
transform in $\om,\bk$ and $\bx$,
\begin{align}
\breve{W}(t,\by,\bq,z) = \int_{\RR}\frac{{\rm d} \om}{2 \pi} \, e^{-i
  \om t} \int_{\RR^d} \frac{{\rm d} \bk}{(2 \pi)^d} \, e^{i \bk \cdot
  \by} \int_{\RR^d} {\rm d} \bx \, e^{- i \bq \cdot \bx}
W(\om,\bk,\bx,z).
\label{eq:Sol1}
\end{align}
Using definition \eqref{eq:TrP2} of the scattering kernel and the
expression \eqref{eq:TR4} of the power spectral density $\hat
\cR_{cc}$, we write
\begin{align*}
Q_{\rm par}(\om,\bk) &= \frac{k_o^2 \sigma_c^2 \ell^{d+1} T}{4}
\int_{\RR} {\rm d}\tilde{t'}\, e^{i \om T \tilde{t'}} \int_{\RR^d}
    {\rm d} \by' \, e^{- i \ell \bk \cdot {\by'}} \int_{\RR} {\rm d}
    z' \, \cR_{cc}(\tilde{t'},{\by'},z') \nonumber \\ &= \frac{k_o^2
      \sigma_c^2 \ell}{4} \int_{\RR} {\rm d}{t}\, e^{i \om {t}}
    \int_{\RR^d} {\rm d} {\by} \, e^{- i \bk \cdot {\by}} \cR \Big(
    \frac{{t}}{T}, \frac{{\by}}{\ell} \Big),
\end{align*}
with $\cR$ defined in \eqref{eq:Sol01}.  Substituting into
\eqref{eq:TRresultP}, we obtain
\begin{align}
\Big[ \partial_z + \frac{\bq}{k_o} \cdot \nabla_\by
  \Big]\breve{W}(t,\by,\bq,z) = \frac{\sigma_c^2 \ell k_o^2}{4}
\Big[ \cR \Big(\frac{t}{T},\frac{\by-\bv t}{\ell} \Big) -
  \cR(0,{\bf 0})\Big] \breve{W}(t,\by,\bq,z),
\label{eq:Sol4}
\end{align}
for $z > 0$, with initial condition obtained from \eqref{eq:paraxini}
and \eqref{eq:Sol1}
\begin{align}
\breve{W}(t,\by,\bq,0) = \breve{W}_0(t,\by) &:= \frac{\sigma_s^2 T_s^2
  \ell_s^{2d}}{4 k_o} \int_{\RR} \frac{{\rm d} \om}{2 \pi} \, e^{-i
  \om t} \int_{\RR^d} \frac{{\rm d} \bk}{(2 \pi)^d} \, e^{i \bk \cdot
  \by} \big|\hat S(T_s \om, \ell_s \bk) \big|^2 \nonumber \\ &=
\frac{\sigma_s^2 T_s\ell_s^{d}}{4 k_o} \int_{\RR} \frac{{\rm d}
  \Omega}{2 \pi} \, e^{-i \Omega \frac{t}{T_s}} \int_{\RR^d}
\frac{{\rm d} \bK}{(2 \pi)^d} \, e^{i \bK \cdot \frac{\by}{\ell_s}}
\big|\hat S(\Omega, \bK) \big|^2.
\label{eq:Sol5}
\end{align}
Note that this condition is independent of $\bq$.

Equation \eqref{eq:Sol4} can be solved by integrating along the
characteristic $ \by = \by_0 + \frac{\bq}{k_o} z, $ stemming from
$\by_0$ at $z = 0$,
\begin{align*}
\breve{W}\Big(t,\by_0 + \frac{\bq}{k_o} z,\bq,z\Big) =&
\breve{W}_0(t,\by_0) \\
  &
  \times \hspace{-0.03in}\exp\hspace{-0.03in}\Big\{\frac{\sigma_c^2 \ell
  k_o^2}{4} \int_0^z \hspace{-0.03in}{\rm d} z' \,
\Big[ \cR\Big(\frac{t}{T},\frac{\by_0 + \frac{\bq}{k_o} z' - \bv_o
    t}{\ell} \Big) - \cR(0,{\bf 0})\Big] \Big\}.
\end{align*}
Substituting $\by_0 = \by - \bq/k_o z$ in this equation, and inverting
the Fourier transform,
\begin{align}
W(\om,\bk,\bx,z) = \int_{\RR} {\rm d} t \, e^{i \om t} \int_{\RR^d}
{\rm d} \by \, e^{-i \bk \cdot \by} \int_{\RR^d} \frac{{\rm d} \bq}{(2
  \pi)^d} \, e^{i \bq \cdot \bx}
\breve{W}_0\Big(t,\by-\frac{\bq}{k_o} z\Big) \nonumber \\ \times
\exp\hspace{-0.03in}\Big\{ \frac{\sigma_c^2 \ell k_o^2}{4} \int_0^z
           {\rm d} z' \, \Big[ \cR\Big(\frac{t}{T},\frac{\by -
               \frac{\bq}{k_o}(z- z') - \bv_o t}{\ell} \Big) -
             \cR(0,{\bf 0})\Big] \Big\}. \label{eq:Sol6}
\end{align}
The result \eqref{eq:Sol0} follows after substituting the expression
\eqref{eq:Sol5} into this equation.

\section{Proof of radiative transfer connection}\label{app:tranproof}  
 
We prove here the result  involving  Eqs. (\ref{eq:T}) and (\ref{eq:rte3d}). 
We start by computing  the different terms in equation \eqref{eq:rte3d}: The first term is 
\[ 
\nabla_{\vec\bk} \Omega(\vec\bk) \cdot \nabla_{\vec\bx}
V(\omega,\vec\bk,\vec\bx) = \frac{c_o}{k_o} \delta\big(
k_z-\beta(\bk)\big)\Big[ \big(\partial_z -\nabla_\bk\beta(\bk)\cdot
  \nabla_\bx\big) W (\omega,\bk,\bx,z) \Big] .
\] 
For the second term we use that if $\vec\bk' = (\bk',\beta(\bk'))$ and
$\vec \bk= (\bk,k_z)$, then
\[ \delta \big( \Omega(\vec\bk)- \Omega(\vec\bk') \big) = \frac{1}{c_o}
\delta \big(|\vec \bk|- k_o\big) = \frac{k_o}{c_o \beta(\bk)}
\delta\big(k_z -\beta(\bk)\big) ,
\]
and 
\begin{align*}
\int_{\RR^{d+1}} \frac{{\rm d} \vec\bk'}{(2\pi)^{d+1}} \int\frac{{\rm
    d}\omega'}{2\pi}
\mathfrak{S}\big(\omega,\omega',\vec\bk,\vec\bk'\big)
V(\omega',\vec\bk',\vec\bx) = \frac{c_o}{k_o} \delta\big(
k_z-\beta(\bk)\big) \nonumber \\ \times \int_{\mathcal{O}}\frac{{\rm
    d} \om' {\rm d} \bk'}{(2 \pi)^{d+1}} Q(\om,\om',\bk,\bk')
W(\omega',\bk',\bx,z) .
\end{align*}
Similarly, for the third term we have that if $\vec\bk =
(\bk,\beta(\bk))$ and $\vec \bk' = (\bk',k'_z)$, then
\[
\delta \big( \Omega(\vec\bk)- \Omega(\vec\bk') \big) = \frac{1}{c_o}
\delta \big(k_o - |\vec \bk| \big) = \frac{k_o}{c_o \beta(\bk')}
\delta\big( k_z' -\beta(\bk') \big) ,
\]
and
\begin{align*}
\int_{\RR^{d+1}} \frac{{\rm d} \vec\bk'}{(2\pi)^{d+1}} \int\frac{{\rm d}\omega'}{2\pi}
\mathfrak{S}\big(\omega,\omega',\vec\bk,\vec\bk'\big) V(\omega,\vec\bk,\vec\bx) =
 \frac{c_o}{k_o} \delta\big(
k_z-\beta(\bk)\big) \nonumber \\ \times \int_{\mathcal{O}}\frac{{\rm d} \om' {\rm d}
  \bk'}{(2 \pi)^{d+1}} Q(\om,\om',\bk,\bk')
  W(\omega,\bk,\bx,z) .
\end{align*}
Gathering the results and using equation \eqref{eq:TRresult}, we
obtain \eqref{eq:rte3d}.

\bibliography{biblio} \bibliographystyle{siam}

\begin{thebibliography}{10}

\bibitem{amzajerdian2013lidar}
{\sc F.~Amzajerdian, D.~Pierrottet, L.~Petway, G.~Hines, V.~Roback, and
  R.~Reisse}, {\em {Lidar sensors for autonomous landing and hazard
  avoidance}}, in Proc. of AIAA Space and Astronautics Forum and Exposition,
  vol.~10, 2013, pp.~6--2013.

\bibitem{BG2016}
{\sc L.~Borcea and J.~Garnier}, {\em Derivation of a one-way radiative transfer
  equation in random media}, Physical Review E, 93 (2016), p.~022115.

\bibitem{BG2016p}
\leavevmode\vrule height 2pt depth -1.6pt width 23pt, {\em Polarization effects
  for electromagnetic wave propagation in random media}, Wave Motion, 63
  (2016), pp.~179--208.

\bibitem{boulanger1995sonic}
{\sc P.~Boulanger, R.~Raspet, and H.~E. Bass}, {\em Sonic boom propagation
  through a realistic turbulent atmosphere}, The Journal of the Acoustical
  Society of America, 98 (1995), pp.~3412--3417.

\bibitem{chandra}
{\sc S.~Chandrasekhar}, {\em Radiative Transfer}, Dover, New York, 1960.

\bibitem{drain1980laser}
{\sc L.~E. Drain}, {\em {The Laser Doppler Techniques}}, Wiley, Chichester,
  1980.

\bibitem{durstprinciples}
{\sc F.~Durst, A.~Melling, and J.~H. Whitelaw}, {\em {Principles and Practice
  of Laser Doppler Anemometry}}, Academic Press, London, 1981.

\bibitem{Mayinger}
{\sc O.~Fieldmann and F.~Mayinger}, eds., {\em {Optical Measurements Techniques
  and Applications}}, Springer-Verlag, Berlin, 2001.

\bibitem{book}
{\sc J.-P. Fouque, J.~Garnier, G.~Papanicolaou, and K.~S{\o}lna}, {\em Wave
  Propagation and Time Reversal in Randomly Layered Media}, Springer, New York,
  2007.

\bibitem{garman2006airborne}
{\sc K.E. Garman, K.~A. Hill, P.~Wyss, M.~Carlsen, J.R. Zimmerman, B.H. Stirm,
  T.Q. Carney, R.~Santini, and P.B. Shepson}, {\em An airborne and wind tunnel
  evaluation of a wind turbulence measurement system for aircraft-based flux
  measurements}, Journal of Atmospheric and Oceanic Technology, 23 (2006),
  pp.~1696--1708.

\bibitem{akira1978wave}
{\sc A.~Ishimaru}, {\em {Wave Propagation and Scattering in Random Media}},
  IEEE Press, Piscataway, 1997.

\bibitem{karyukin1982influence}
{\sc G.~A. Karyukin}, {\em Influence of wind on operation of radar acoustic
  atmospheric sounding systems}, Izv. Acad. Sci. USSR, Atmos. Oceanic Phys, 18
  (1982), pp.~26--30.

\bibitem{molland2017ship}
{\sc A.~F. Molland, S.~R. Turnock, and D.~A. Hudson}, {\em Ship Resistance and
  Propulsion}, Cambridge University Press, Cambridge, 2017.

\bibitem{munk2009ocean}
{\sc W.~Munk, P.~Worcester, and C.~Wunsch}, {\em Ocean Acoustic Tomography},
  Cambridge University Press, Cambridge, 2009.

\bibitem{oberg1990laser}
{\sc P.~A. Oberg}, {\em Laser-doppler flowmetry}, Critical Reviews in
  Biomedical Engineering, 18 (1990), pp.~125--163.

\bibitem{ostashev2015acoustics}
{\sc V.~E. Ostashev and D.~K. Wilson}, {\em {Acoustics in moving inhomogeneous
  media}}, CRC Press, Boca Raton, 2015.

\bibitem{PK1974}
{\sc G.~Papanicolaou and W.~Kohler}, {\em {Asymptotic theory of mixing
  stochastic ordinary differential equations}}, Communications on Pure and
  Applied Mathematics, 27 (1974), pp.~641--668.

\bibitem{PW1994}
{\sc G.~Papanicolaou and S.~Weinryb}, {\em A functional limit theorem for waves
  reflected by a random medium}, Applied Mathematics and Optimization, 30
  (1994), pp.~307--334.

\bibitem{pierce1990wave}
{\sc A.~D. Pierce}, {\em {Wave equation for sound in fluids with unsteady
  inhomogeneous flow}}, {The Journal of the Acoustical Society of America}, 87
  (1990), pp.~2292--2299.

\bibitem{ryzhik96}
{\sc L.~Ryzhik, G.~Papanicolaou, and J.B. Keller}, {\em Transport equations for
  elastic and other waves in random media}, Wave motion, 24 (1996),
  pp.~327--370.

\bibitem{schomer1983noise}
{\sc P.~D. Schomer}, {\em Noise monitoring in the vicinity of general aviation
  airports}, The Journal of the Acoustical Society of America, 74 (1983),
  pp.~1764--1772.

\end{thebibliography}
\end{document}